\documentclass[12pt]{amsart}
\usepackage{amscd,amssymb,graphics,color,a4wide,hyperref}
\usepackage{graphicx}
\usepackage{psfrag} 
\usepackage{marvosym}

\usepackage{mathrsfs}
\input xy
\xyoption{all}

\DeclareFontFamily{U}{rsf}{}
\DeclareFontShape{U}{rsf}{m}{n}{
  <5> <6> rsfs5 <7> <8> <9> rsfs7 <10->  rsfs10}{}
\DeclareMathAlphabet{\mathscr}{U}{rsf}{m}{n}

\newtheorem{theorem}{Theorem}[section]

\newtheorem{proposition}[theorem]{Proposition}

\newtheorem{conjecture}[theorem]{Conjecture}

\theoremstyle{definition}
\newtheorem{definition}[theorem]{Definition}

\newtheorem{example}[theorem]{Example}
\newtheorem{examples}[theorem]{Examples}

\theoremstyle{remark}

\numberwithin{equation}{section}

\newcommand{\NN} {\mathbb{N}}
\newcommand{\ZZ} {\mathbb{Z}}

\newcommand{\RR} {\mathbb{R}}
\newcommand{\VV} {\mathbb{V}}
\newcommand{\CC} {\mathbb{C}}

\newcommand{\PP} {\mathbb{P}}
\renewcommand{\AA} {\mathbb{A}}
\newcommand{\GG} {\mathbb{G}}

\newcommand {\shA} {\mathcal{A}}
\newcommand {\shAff} {\mathcal{A}\text{\textit{ff}}}

\newcommand {\shL} {\mathcal{L}}

\newcommand {\shT} {\mathcal{T}}

\newcommand {\shP} {\mathcal{P}}

\newcommand {\shX} {\mathcal{X}}
\newcommand {\shY} {\mathcal{Y}}

\newcommand {\foD} {\mathfrak{D}}
\newcommand {\fod} {\mathfrak{d}}

\newcommand {\fom} {\mathfrak{m}}

\newcommand {\fop} {\mathfrak{p}}

\newcommand {\foX} {\mathfrak{X}}

\newcommand {\Aff} {\operatorname{Aff}}

\newcommand {\Aut} {\operatorname{Aut}}

\newcommand {\dual} {\vee}

\newcommand {\ev} {\mathrm{ev}}

\newcommand {\Fuk} {\mathrm{Fuk}}

\newcommand {\GL} {\operatorname{GL}}

\newcommand {\Hom} {\operatorname{Hom}}

\renewcommand {\Im} {\operatorname{Im}}

\newcommand {\Int} {\operatorname{Int}}

\renewcommand {\ker } {\operatorname{ker}}

\newcommand {\Lift} {\operatorname{Lift}}

\newcommand {\lra} {\longrightarrow}

\newcommand {\Mono} {\operatorname{Mono}}

\renewcommand{\O} {\mathcal{O}}

\newcommand {\ord} {\operatorname{ord}}
\newcommand {\out} {\mathrm{out}}

\renewcommand{\P} {\mathscr{P}}

\newcommand {\Sing} {\operatorname{Sing}}
\newcommand {\SL} {\operatorname{SL}}
\newcommand {\Spec} {\operatorname{Spec}}
\newcommand {\Spf} {\operatorname{Spf}}

\newcommand {\trcone} {\overline{C}}

\newcommand {\ul} {\underline}

\newcommand{\bbfamily}{\fontencoding{U}\fontfamily{bbold}\selectfont}
\newcommand{\textbb}[1]{{\bbfamily#1}}
\newcommand {\lfor} {\mbox{\textbb{[}}}
\newcommand {\rfor} {\mbox{\textbb{]}}}
\newcommand {\T} {\shT}
\newcommand {\X} {\shX}

\newcommand {\Gm} {\GG_m}

\def\mydate{\ifcase\month \or January\or February\or March\or
April\or May\or June\or July\or August\or September\or October\or 
November\or December\fi \space\number\day,\space\number\year}

\newlength{\picwidth} 
\newlength{\miniwidth} 
\newlength{\nanowidth} 
\newlength{\melowidth} 
\newlength{\leftminiwidth} 
\newlength{\rightminiwidth} 
\newlength{\minipagewidth} 

\begin{document}
\def\mapright#1{\smash{
 \mathop{\longrightarrow}\limits^{#1}}}
\def\mapleft#1{\smash{
 \mathop{\longleftarrow}\limits^{#1}}}
\def\exact#1#2#3{0\to#1\to#2\to#3\to0}
\def\mapup#1{\Big\uparrow
  \rlap{$\vcenter{\hbox{$\scriptstyle#1$}}$}}
\def\mapdown#1{\Big\downarrow
  \rlap{$\vcenter{\hbox{$\scriptstyle#1$}}$}}
\def\dual#1{{#1}^{\scriptscriptstyle \vee}}
\def\invlim{\mathop{\rm lim}\limits_{\longleftarrow}}
\def\rto{\raise.5ex\hbox{$\scriptscriptstyle ---\!\!\!>$}}

\input epsf.tex
\title
[Theta functions and mirror symmetry]
{Theta functions and mirror symmetry}
\author{Mark Gross}
\address{UCSD Mathematics, 9500 Gilman Drive, La Jolla, CA 92093-0112, USA}
\email{mgross@math.ucsd.edu}
\thanks{This work was partially supported by NSF grants
0854987 and 1105871.}

\author{Bernd Siebert}
\address{Department Mathematik,
Universit\"at Hamburg, Bundesstra\ss e~55, 20146 Hamburg,
Germany}
\email{bernd.siebert@math.uni-hamburg.de}

\maketitle
\tableofcontents
\bigskip


\section*{Introduction}

The classical subject of theta functions has a very rich history
dating back to the nineteenth century. In modern algebraic
geometry, they arise as sections of ample line bundles
on abelian varieties, canonically defined after making
some discrete choices of data. The definition of theta functions
depends fundamentally on the group law,
leaving the impression that they are a feature restricted to
abelian varieties. However, new insights from mirror symmetry suggest
that they exist much more generally, even on some of 
the most familiar varieties.

Mirror symmetry began as a phenomenon in string theory in
1989, with the suggestion that Calabi-Yau manifolds should
come in pairs. Work of Greene and Plesser \cite{GrPl} 
and Candelas, Lynker and Schimmrigk \cite{CLS} gave the first hint that
there was mathematical justification for this idea, with constructions
given of pairs of Calabi-Yau three-folds $X$, $\check X$, with the property
that $\chi(X)=-\chi(\check X)$. More precisely, the Hodge numbers
of these pairs obey the relation
\[
h^{1,1}(X)=h^{1,2}(\check X),\quad h^{1,2}(X)=h^{1,1}(\check X).
\]
In 1991, Candelas, de la Ossa, Green and Parkes \cite{COGP} achieved an
astonishing breakthrough in exploring the mathematical ramifications
of some string-theoretic predictions. In particular, using string theory
as a guideline, they carried out certain period integral calculations
for the mirror of the quintic three-fold in $\CC\PP^4$, and obtained
a generating function for the numbers $N_d$, $d\ge 1$, where $N_d$
is the number of rational curves of degree $d$ in the quintic three-fold.

This immediately attracted attention from mathematicians, and there has
followed twenty years of very rewarding efforts to understand the mathematics
underlying mirror symmetry. 

A great deal of progress has been made, but much remains to be done. 
In this survey article, we will discuss certain aspects of this search
for understanding, guided by a relationship between mirror symmetry
and theta functions. In particular, we will discuss the surprising
implication of mirror symmetry that near large complex structure limits
in complex moduli space, Calabi-Yau manifolds also carry
theta functions. Roughly speaking, these will be
canonically defined bases for the space of sections of line bundles. 
This implication was first suggested, as far as we know, by the
late Andrei Tyurin, see \cite{Ty99}.

The constructions in fact apply much more broadly than just to Calabi-Yau
manifolds. For example, in \cite{GHKII} (see Sean Keel's lecture
\cite{Ke11}), it is proven that some of the most familiar
varieties in algebraic geometry, including many familiar affine
rational surfaces, carry theta functions. 
In addition, much stronger results apply to the surface
case, with \cite{GHKK3} proving a strong form of Tyurin's conjecture.
However, in this survey we will focus only on the simplest aspects of the
construction.

The discussion here represents a distillation of a number of ongoing
joint projects with varying groups of coauthors. The 
relationship between integral affine manifolds,
degenerations of Calabi-Yau varieties and mirror symmetry discussed
here is based on a long-term project of the authors of this survey.
The construction of theta functions as described here has come out
of joint work with Paul Hacking and Sean Keel,
while the relationship between theta functions and homological mirror
symmetry is based on forthcoming work of Mohammed Abouzaid with Gross
and Siebert.

This survey is based on a lecture delivered by the first author
at the JDG 2011 conference in April 2011 at Harvard University.
We would like to thank Professor Yau for this invitation and the
opportunity to contribute to the proceedings. We also thank our coworkers
on various projects described here: Mohammed Abouzaid, Paul Hacking and
Sean Keel.

\section{The geometry of mirror symmetry: HMS and SYZ}

There are two principal approaches to the geometry underlying
mirror symmetry: Kontsevich's homological mirror symmetry conjecture
(HMS) \cite{K95} and the Strominger-Yau-Zaslow (SYZ) conjecture \cite{SYZ}.
Taken together, they suggest the existence of theta functions.

These conjectures are as follows. Consider a mirror pair of Calabi-Yau
manifolds, $X$ and $\check X$. To be somewhat more
precise, we should consider
Calabi-Yau manifolds with Ricci-flat K\"ahler metric, so that mirror
symmetry is an involution
\[
(X,J,\omega)\leftrightarrow (\check X,\check J,\check\omega).
\]
Here $J$ is the complex structure and $\omega$ the K\"ahler 
form on $X$. One expects that $J$
determines the K\"ahler structure $\check\omega$ and $\omega$
determines the complex structure $\check J$.\footnote{This discussion
ignores the $B$-field.} Kontsevich's fundamental insight is that
the isomorphism that mirror symmetry predicts between the complex geometry
of $(X,J)$ and the symplectic geometry of $(\check X,\check\omega)$ can
be expressed in a categorical setting:

\begin{conjecture}[Homological mirror symmetry] There is an equivalence
of categories between the derived category $D^b(X)$ of bounded complexes
of coherent sheaves on $X$ and $\Fuk(\check X)$, the Fukaya category
of Lagrangian submanifolds on $\check X$.
\end{conjecture}

There are many technical issues hiding in this statement, the least
of which is showing that the Fukaya category makes sense. In particular,
$\Fuk(\check X)$ is not a category in the traditional sense, as composition
of morphisms is not associative. Rather, it is an $A_{\infty}$-category,
which essentially means that there is a sequence of higher composition
maps which measure the failure of associativity; we will be more precise
shortly.

To first approximation, the objects of $\Fuk(\check X)$ are Lagrangian
submanifolds of $\check X$, i.e., submanifolds $L\subseteq \check X$
with $\dim_{\RR} L=\dim_{\CC} \check X$ and $\check\omega|_L=0$. We
define the $\Hom$ between objects as follows. Let $\Lambda$ be
the Novikov ring, i.e., the ring of power series $\sum_{i=1}^{\infty}
a_iq^{r_i}$ where $a_i\in\CC$, $r_i\in \RR_{\ge 0}$, $r_i\rightarrow\infty$
as $i\rightarrow\infty$.
Given two Lagrangian submanifolds $L_0,L_1$, we define
\[
\Hom(L_0,L_1)=\bigoplus_{p\in L_0\cap L_1} \Lambda [p],
\]
assuming that $L_0$ and $L_1$ intersect transversally (if not, we
can perturb one of them via a generic Hamiltonian isotopy). In fact,
this is a graded $\Lambda$-module, with the degree of $[p]$ being
the so-called \emph{Maslov index} of $p$. One can then define a series
of maps
\[
\mu_d:\Hom(L_{d-1},L_d)\otimes\cdots\otimes \Hom(L_0,L_1)\rightarrow
\Hom(L_0,L_d)
\]
for $d\ge 1$, $L_0,\ldots,L_d$ Lagrangian submanifolds of $\check X$.
Roughly this map is defined by counting certain holomorphic disks:
\[
\mu_d(p_{d-1,d},\ldots,p_{0,1})
=\sum_{p_{0,d}\in L_0\cap L_d}\sum_{\psi:D\rightarrow \check X} 
\pm q^{\int_D \psi^*\check\omega} [p_{0,d}]
\]
where the second sum is over all holomorphic maps $\psi:D\rightarrow
X$ such that there are cyclically ordered points $t_0,\ldots,t_d\in\partial
D$ with $\psi(t_i)=p_{i,i+1}$, $\psi(t_d)=p_{0,d}$, and
$\psi([t_i,t_{i+1}])\subseteq L_{i+1}$ and $\psi([t_d,t_0])
\subseteq L_0$. Here $[t_i,t_{i+1}]$ denotes the interval on $\partial D$
between $t_i$ and $t_{i+1}$. See Figure \ref{Ainfinitydisk}. The contribution to
$\mu_d(p_{d-1,d},\ldots,p_{0,1})$ from $p_{0,d}$ is only counted if the
expected dimension of the moduli space of disks is zero; this makes
$\mu_d$ into a chain map of degree $2-d$. In suitably nice cases, 
these operations will satisfy the so-called $A_{\infty}$-relations,
which are
\[
\sum_{1\le p\le d\atop 0\le q\le d-p}
\pm \mu_{d-p+1}(a_d,\ldots,a_{p+q+1},\mu_p(a_{p+q},\ldots,a_{q+1}),
a_q,\ldots,a_1)=0.
\]
This tells us that $\mu_1$ turns $\Hom(L_0,L_1)$ into a complex,
that $\mu_2$ is associative up to homotopy, and so on. Since $\mu_2$
will play the role of composition of morphisms for us, this means
the Fukaya category is not in general a category, because composition
is not associative.
\begin{figure}
\input{Ainfinitydisk.pstex_t}
\caption{}
\label{Ainfinitydisk}
\end{figure}

One might now object that $D^b(X)$ is a genuine category, so in what
sense is it isomorphic to something which is not a genuine category?
It turns out that there is a natural way to put an $A_{\infty}$-category
structure on $D^b(X)$ which is really enriching the structure of
$D^b(X)$; the statement of HMS then says that we expect $D^b(X)$
and $\Fuk(\check X)$ to be quasi-isomorphic as $A_{\infty}$-categories,
which is a well-defined notion.

We will not go into more technical details about HMS, as we shall only
need a small part of it. Instead, we move on to discuss the 
Strominger-Yau-Zaslow (SYZ) conjecture \cite{SYZ}, dating from 1996. This
was the first idea giving a truly geometric interpretation of
mirror symmetry. Fix now an $n$-dimensional Calabi-Yau manifold $X$ with a
nowhere vanishing holomorphic $n$-form $\Omega$ and a symplectic form
$\omega$: we say a submanifold $M\subseteq X$ is \emph{special Lagrangian}
if $M$ is Lagrangian with respect to $\omega$ and furthermore $\Im\Omega|_M=0$.
This notion was introduced by Harvey and Lawson in \cite{HL82}; special
Lagrangian submanifolds are volume minimizing in their homology class.
Suppose $X$ has a mirror $\check X$, with holomorphic $n$-form $\check\Omega$
and symplectic form $\check\omega$.
We then have:

\begin{conjecture}[The Strominger-Yau-Zaslow conjecture]
There are continuous maps $f:X\rightarrow B$, $\check f:\check X\rightarrow
B$ whose fibres are special Lagrangian, and whose general fibres
are dual $n$-tori.
\end{conjecture}

This is a purposefully vague statement, partly because we are very far
from a proof of anything resembling this conjecture: see 
\cite{Gr98},\cite{Gr99},\cite{Gr00} for 
detailed discussion of more precise forms of this conjecture.
Let us just say at this point that the duality implies that the
topological monodromy of the smooth part of $f$ is the transpose of
the topological monodromy of $\check f$.

Early work aimed at understanding the conjecture includes
\cite{Gr98},\cite{Gr99},\cite{Hi97}.
In particular in \cite{Gr98}, the first author conjectured
that Lagrangian sections of $\check f$ should be expected to be mirror,
under HMS, to line bundles on $X$. A more precise correspondence
was predicted there, with specific predictions on which topological
isotopy class of sections corresponded to which numerical equivalence
class of line bundles. This idea was used in a number of different
situations: for example, work of Polishchuk and Zaslow \cite{PZ98}
give an explicit
correspondence between special Lagrangian sections on the obvious
SYZ fibration on an elliptic curve and line bundles on the mirror
elliptic curve. We shall say more about this in \S \ref{abeliansection}.

The fundamental idea we shall pursue in this paper is the following.
Suppose $L_0$ is a Lagrangian section of $\check f$ corresponding to 
the structure sheaf $\O_X$, and $L_1$ is a Lagrangian section
corresponding to an ample line bundle $\shL$ on $X$. Then HMS
should yield an isomorphism $\Hom(L_0,L_1)\cong \Hom(\O_X,\shL)\otimes_{\CC}
\Lambda$. Here the $\Hom$'s are in the Fukaya and derived categories
respectively, but after taking cohomology, one expects on the right
to only get a contribution from $H^0(X,\shL)$ as all higher cohomology
vanishes. If all the intersection points of $L_0$ and $L_1$
are Maslov index zero, i.e., if somehow the intersection is particularly nice
so that there is no $\mu_1$, then one has of course a basis for
$\Hom(L_0,L_1)$ given by these intersection points, and these
correspond to elements of $H^0(X,\shL)$.

The moral of this is: suppose we have particularly canonical choices
of Lagrangian sections corresponding to ample line bundles. Then HMS
predicts the existence of a canonical basis of sections of ample
line bundles. We are going to call elements of such a canonical
basis \emph{theta functions}.

To the best of our knowledge, the existence of such a canonical basis
was first suggested by the late Andrei Tyurin; the first author
heard him speak about these in a lectures series at the University
of Warwick in 1999: see especially the remark on p.\,36 of \cite{Ty99}.

We will first make more precise
what these canonical sections should be, and then in the next
section argue in the case of abelian varieties that the corresponding
basis indeed coincides with classical theta functions. Before doing so,
we need to explain more structure underlying the correct way of thinking
about the SYZ conjecture.

What is actually most important about the SYZ conjecture are certain
structures which should appear on the base $B$ of the two dual
fibrations. Suppose we have fibrations as in the conjecture.
Let $\Delta\subseteq B$ be the set of critical values
of $f$ and $B_0:=B\setminus \Delta$, so that if $x\in B_0$, 
$f^{-1}(x)$ is a smooth torus. It was first observed by Hitchin in \cite{Hi97} 
that the special Lagrangian fibration induces two different affine structures
on $B_0$, one induced by $\omega$ (this affine structure arises
from the Arnold-Liouville theorem)
and one induced by $\Im\Omega$. A precise definition:

\begin{definition}
An affine structure on an $n$-dimensional 
real manifold $B$ is a set of coordinate
charts $\{\psi_i:U_i\rightarrow \RR^n\}$ on an open cover $\{U_i\}$
of $B$ whose transition maps $\psi_j\circ\psi_i^{-1}$ lie in $\Aff(\RR^n)$,
the affine linear group of $\RR^n$. We say the structure is \emph{tropical}
if the transition maps lie in $\RR^n\rtimes \GL_n(\ZZ)$, and
\emph{integral} if the transition maps lie in $\Aff(\ZZ^n)$.

A \emph{(tropical, integral) affine manifold with singularities} is
a manifold $B$ along with a codimension $\ge 2$ subset $\Delta\subseteq B$
and a (tropical, integral) affine structure on $B_0:=B\setminus \Delta$.
\end{definition}

In fact, the affine structures induced by special Lagrangian fibrations
are tropical, so we obtain tropical affine manifolds with singularities
(except for the fact that genuine special Lagrangian fibrations are
expected to have codimension one discriminant loci which retract
onto a codimension two subset, as demonstrated
by Joyce in many examples \cite{J03}). 

It is convenient now to largely forget about special Lagrangian fibrations,
as we don't know if they exist, and instead focus on the tropical affine
manifolds arising from them. It is in fact now fairly well understood
what such manifolds should look like, even if the fibrations aren't known!
See for example \cite{Gr09}.

In fact, tropical affine manifolds quickly give rise to a toy version
of mirror symmetry:

\begin{definition}
If $B$ is a tropical affine manifold, then let $\Lambda$
be the local system contained in the tangent bundle $\shT_B$
given locally by integral linear combinations
of coordinate vector fields $\partial/\partial y_1,\ldots,\partial/
\partial y_n$, where $y_1,\ldots,y_n$ are local tropical affine coordinates.
The fact that transition maps lie in $\RR^n\rtimes \GL_n(\ZZ)$ rather
than $\Aff(\RR^n)$ says this local system is well-defined, independently
of coordinates. Similarly, let $\check\Lambda\subseteq \T_B^*$
be the local system given locally by integral linear combinations
of $dy_1,\ldots,dy_n$. Set
\begin{align*}
X(B):= {} & \T_B/\Lambda,\\
\check X(B):= {} & \T^*_B/\check\Lambda.
\end{align*}
We have projections $f:X(B)\rightarrow B$ and $\check f:\check X(B)\rightarrow
B$ which are dual torus fibrations. 
\end{definition}

Note that $X(B)$ comes along with a natural complex structure. This
is most easily described by specifying the almost complex structure $J$.
There is a natural flat connection on $\T_B$ such that sections of
$\Lambda$ are flat sections. At any point in $\T_B$, the horizontal
and vertical tangent spaces are both isomorphic to the tangent space
to $B$, and $J$ interchanges these two spaces, inserting an appropriate
sign to ensure $J^2=-1$. It is easy to see that this structure is
integrable, identifying $f^{-1}(U)$, for $U\subseteq B$ a small open set,
with a $T^n$-invariant open subset of $(\CC^*)^n$. 

Furthermore, $\check X(B)$ carries a natural symplectic structure:
as always, $\T_B^*$ carries a canonical symplectic form, and one checks
it descends to the quotient.

As a consequence, we can view the correspondence $X(B)\leftrightarrow
\check X(B)$ as a toy version of mirror symmetry. In this discussion
we see half of mirror symmetry, as we don't have a symplectic structure
on $X(B)$ or a complex structure on $\check X(B)$. 

How close is this correspondence to actual mirror symmetry? If $B$
is compact, e.g., $B=\RR^n/\Gamma$ for a lattice $\Gamma$, then $X(B)$
is a complex torus, and the toy description gives a completely
satisfactory description of mirror symmetry; we shall make use of this
in \S \ref{abeliansection}. 
However, in general, one should work with $B$ a tropical affine
manifold with singularities, in which case one only has a subset $B_0\subseteq
B$ with an affine structure. So one can then ask to what extent can
one compactify $X(B_0)$ or $\check X(B_0)$. One has the following
general observations:
\begin{enumerate}
\item In various nice cases, $X(B_0)$ and $\check X(B_0)$ can be compactified
topologically: see \cite{Gr01} 
for the three-dimensional case, and work in progress
\cite{GS13} for similar results in all dimensions. In particular, \cite{Gr01}
gives a complete description of the quintic threefold and its mirror
from this point of view.
\item In various nice cases in dimensions two and three, $\check X(B_0)$
can be compactified to some $\check X(B)$ 
in the symplectic category, see \cite{CBM09}.
\item As a complex manifold, $X(B_0)$ can almost never be compactified.
This is a crucial point for mirror symmetry. There are \emph{instanton
corrections} that one needs to make to the complex structure on
$X(B_0)$ before one can hope to compactify this. This was first
explored by Fukaya in \cite{F05}, in the two-dimensional (K3) case. 
That paper was the first to suggest the philosophy: \emph{the corrections
to the complex structure on $X(B_0)$ arise from pseudo-holomorphic disks
in $\check X(B)$ with
boundary on fibres of the SYZ fibration}.
\end{enumerate}

The direct analytic approach of Fukaya suffers from huge technical
difficulties and as such was only a heuristic. 
To carry out this program in a more practical way, 
a switch to an easier category is necessary. 
Kontsevich and Soibelman used the rigid analytic category, constructing
in \cite{KS06} a rigid analytic K3 surface from a tropical affine surface
with $24$ singular points. In parallel, we had been working on a program
to address this problem in all dimensions using logarithmic geometry;
combining our approach with some ideas of \cite{KS06}, we provided a solution
to this problem in all dimensions, in a somewhat different category.
Roughly, our result shown in \cite{GS11} is as follows.

\begin{theorem}
\label{mainsmoothingtheorem}
Suppose given an integral affine manifold with singularities
$B$. Suppose furthermore that the singularities are ``nice'' and $B$ comes
with a decomposition $\P$ into lattice polytopes. Suppose furthermore
given a strictly convex, multi-valued piecewise linear function $\varphi$
on $B$. Then one can construct a one-parameter flat family 
$\pi:\check\shX\rightarrow
\Spec\CC\lfor t\rfor$ from this data whose central fibre is
\[
\check\shX_0=\bigcup_{\sigma\in\P_{\max}} \PP_{\sigma},
\]
where $\P_{\max}$ is the set of maximal cells in the polyhedral decomposition,
and $\PP_{\sigma}$ is the projective toric variety defined by the lattice
polytope $\sigma$. These toric varieties are glued together along toric
strata as dictated by the combinatorics of $\P$.
Furthermore, $\check\shX$ comes along with a relatively ample
line bundle $\shL$.
\end{theorem}

There are a number of important features of this construction:
\begin{enumerate}
\item It is an explicit construction, giving an order-by-order algorithm
for gluing standard thickenings of affine pieces of the irreducible
components of $\check\shX_0$. This data is described by what we call a 
\emph{structure}, and as we shall see, is really controlled by
counts of holomorphic disks on the mirror side.
\item
There is
a notion of discrete Legendre transform which allows one to associate
to the triple of data $(B,\P,\varphi)$ another triple $(\check B,
\check\P,\check\varphi)$. The polyhedral decomposition $\check\P$ is
dual to $\P$, and the affine structure on $\check B$ is dual to that
on $B$ in some precise sense, see \cite{GS06}, \S1.4. 
Then applying the above theorem
to the dual data yields the mirror Calabi-Yau.
\item The family constructed in Theorem \ref{mainsmoothingtheorem} 
can be extended to a flat family of complex
analytic spaces ${\mathscr X}$ over a disk $D$. A general fibre
of this family, ${\mathscr X}_t$, has a K\"ahler form repesenting
$c_1(\shL)$. The expectation is that this symplectic manifold
is a compactification of $\check X(B_0)$. Furthermore, as a complex
manifold, it should roughly be a compactification of a small deformation
of the complex structure on $X_{\epsilon}(\check B_0)$,
where $X_{\epsilon}(\check B_0)=\T_{\check 
B_0}/\epsilon\Lambda$ and $\epsilon>0$
is a real number. See \cite{GS03} for some details of this; more
details will appear in \cite{GS13}. 
\end{enumerate}

There is one confusing point in this discussion: 
the role of $X(B)$ and $\check X(B)$
has been interchanged. We originally said we wanted to compactify
$X(B_0)$. Instead, we compactified $\check X(B_0)$. This makes
sense from several points of view. 

First, we have constructed
a whole family of complex manifolds, but they are symplectomorphic
as symplectic manifolds. So it makes sense that we get $\check X(B)$,
which comes with a canonical symplectic structure. The integrality
of the affine structure on $B$ guarantees that the symplectic form
on $\check X(B)$ represents an integral cohomology class.

Second, if one doesn't like this switch, then one can work with
the Legendre dual manifold $\check B$. The work of Fukaya \cite{F05} 
and Kontsevich and Soibelman \cite{KS06} did precisely this. But it turns
out that structures are nicer objects on $B$ than on $\check B$.
On $B$, the data controlling the family $\check\shX$, the structure,
is essentially tropical in nature, and can be viewed as a union of
tropical trees on $B$. If one works on $\check B$, one instead needs to
use trees made of gradient flow lines, and this can produce some technical
difficulties. Working on $B$ makes many aspects of our work effective.

\medskip

Let us now return to theta functions. We note our construction comes
with a canonical ample line bundle, whose first Chern class is
represented by the symplectic form on $\check X(B)$. Now a line bundle
should be mirror to a section of the SYZ fibration, so it is natural
to ask whether $X(B)\rightarrow B$ comes with a natural section. Since we
haven't given an explicit description of the compactification in this
paper, let us at least
answer this question over $B_0$. 
There is in fact a whole set of natural sections, indexed by $\ell\in \ZZ$, given
in local integral affine coordinates by
\begin{equation}
\label{sectiondef}
\sigma_\ell:
(y_1,\ldots,y_n)\mapsto -\sum_{i=1}^n \ell\cdot y_i{\partial\over\partial y_i}.
\end{equation}
Note that modulo integral vector fields, i.e., sections of $\Lambda$,
this vector field is well-defined independently of the choice of
integral affine coordinates. Call the image of this section $L_\ell$.

Of course there is no symplectic structure on $X(B_0)$, so it doesn't
quite make sense to call these Lagrangian sections, but one can imagine
that one can find symplectic structures which make these sections Lagrangian,
and then deduce some consequences. 

The most important consequence 
is the description of the set $L_0\cap L_\ell$, namely
\[
f(L_0\cap L_\ell)=B_0\left({1\over \ell}\ZZ\right),
\]
where the latter denotes the set of points of $B_0$ whose coordinates
in any (hence all) integral affine coordinate charts lie in ${1\over \ell}\ZZ$.

Let us hypothesize that $L_\ell$ is mirror to the line bundle $\shL^\ell$.
Since for an ample line bundle on a Calabi-Yau manifold all higher cohomology
vanishes, we are led to the following conjecture:

\begin{conjecture} 
For $\ell>0$, there is a canonical
basis of $\Gamma(\check\shX,\shL^\ell)$ as a $\CC\lfor t\rfor$-module
indexed by elements of $B({1\over \ell}\ZZ)$.
\end{conjecture}

In the sections that follow, we will first explain why theta functions
for abelian varieties fit naturally into such a conjecture. Next, we 
outline the proof of this conjecture given in \cite{GHKS}, with
the precise statement given in Theorem \ref{GHKSmain}. We finally 
explain various applications of the existence of such a basis.

\section{Theta functions for abelian varieties and the Mumford
construction}
\label{abeliansection}

In the case that $B$ is a torus, our construction in fact recovers
Mumford's 
description of degenerations of abelian varieties \cite{M72}, see also
\cite{AN99}. Theorem \ref{mainsmoothingtheorem} can be viewed as a vast
generalization of this construction. 
We will briefly review a simple version of Mumford's construction. 

The starting data is a lattice $M\cong\ZZ^n$, $M_{\RR}=M\otimes_{\ZZ}\RR$,
$N=\Hom_{\ZZ}(M,\ZZ)$,
a sublattice $\Gamma\subseteq M$, a $\Gamma$-periodic polyhedral
decomposition $\P$ of $M_{\RR}$, and a strictly convex
piecewise linear function with integral slopes $\varphi:M_{\RR}\rightarrow
\RR$ satisfying a periodicity condition, for $\gamma\in\Gamma$,
\[
\varphi(m+\gamma)=\varphi(m)+\alpha_{\gamma}(m)
\]
for some affine linear function $\alpha_{\gamma}$ depending on $\gamma$.
The affine manifold $B$ in Theorem \ref{mainsmoothingtheorem}
will be $M_{\RR}/\Gamma$ in this setup.

From this data one builds an unbounded polyhedron in $M_{\RR}\oplus\RR$:
\[
\Delta_{\varphi}:=\{(m,r)\,|\, m\in M_{\RR}, r\ge \varphi(m)\}.
\]
The normal fan of this polyhedron in $N_{\RR}\oplus\RR$ 
is a fan $\Sigma_{\varphi}$ with an infinite number of cones, defining a 
toric variety $X_{\varphi}$ which is not of finite type. Note that the
one-dimensional rays of $\Sigma_{\varphi}$ are in one-to-one correspondence
with the maximal cells $\sigma$ of $\P$; if $n_{\sigma}\in N$ is the slope
of $\varphi|_{\sigma}$, then $(-n_{\sigma},1)$ is the corresponding ray
in $\Sigma_{\varphi}$. Further, $\Gamma$ acts on $N\oplus\ZZ$; indeed,
$\gamma\in
\Gamma$ acts by taking $(n,r)\mapsto (n-r\cdot d\alpha_{\gamma},r)$,
where $d\alpha_{\gamma}$ denotes the differential, or, 
slope, of $\alpha_{\gamma}$. This action preserves $\Sigma_{\varphi}$. 

The projection $N_{\RR}\oplus\RR\rightarrow \RR$ defines a map
$\pi:X_{\varphi}\rightarrow \AA^1$. The fibres of this map are
algebraic tori $(\CC^*)^n$ except for $\pi^{-1}(0)$, which is an infinite
union of proper toric varieties. Furthermore, the action of $\Gamma$
preserves this map, and yields an action of $\Gamma$ on the irreducible
components of $\pi^{-1}(0)$.

While the $\Gamma$-action is global,
it does not act properly discontinuously except on the subset $\pi^{-1}(D)$,
where $D\subseteq\AA^1$ is the unit disk. Thus we get a family
\[ 
\pi:\pi^{-1}(D)/\Gamma\rightarrow D
\]
whose general fibre is an abelian variety and such that the fibre
over zero is a union of toric varieties. 

We would actually prefer to work formally here, and instead consider
\[
\shA:=(X_{\Sigma}\times_{\AA^1}\Spec\CC\lfor t\rfor)/\Gamma.
\]
The quotient can be taken by dividing out the 
formal completion of $X_{\varphi}$
along $\pi^{-1}(0)$ by the action of $\Gamma$, then showing that there
is an ample line bundle on this quotient, and finally applying
Grothendieck existence to get a scheme over $\Spec\CC\lfor t\rfor$. In
fact, the existence of the ample line bundle will follow from the discussion
below.

The family $\shA\rightarrow\Spec\CC\lfor t\rfor$
is precisely the family produced by
Theorem \ref{mainsmoothingtheorem} from the data
$B=M_{\RR}/\Gamma$, polyhedral decomposition given by the image of
$\P$ in $B$, and multi-valued piecewise linear function $\varphi$
as given.

We now would like to understand traditional theta functions in this
context. As already understood in \cite{M72}, one observes
that the polyhedron $\Delta_{\varphi}$ induces a line bundle $\shL$
on $X_{\Sigma}$, and the bundle $\shL$ descends to the quotients
\[
\shA_k:=(X_{\Sigma}\times_{\AA^1}\Spec\CC[t]/(t^{k+1}))/\Gamma
\]
of the $k$-th order thickenings of the central fibre of $X_{\Sigma}\rightarrow
\AA^1$. To show $\shL$ descends, one just needs to define an integral
linear action
of $\Gamma$ on the cone $C(\Delta_{\varphi})\subseteq M_{\RR}\oplus\RR\oplus\RR$
defined as
\[
C(\Delta_{\varphi})=\overline{\{(\ell m,\ell r,\ell)\,|\,(m,r)\in\Delta_{\varphi}, 
\ell\in \RR_{\ge 0}\}}.
\]
Note taking the closure just adds $\{0\}\times\RR\times\{0\}$ to the set.
If $c_{\gamma}$ is the constant part of $\alpha_{\gamma}$ and
$d\alpha_{\gamma}$ the differential of $\alpha_{\gamma}$, (or equivalently,
$d\alpha_{\gamma}$ is
the linear part of $\alpha_{\gamma}$), one checks such an action is given by
$\gamma \mapsto \psi_{\gamma}\in \Aut(M\oplus\ZZ
\oplus\ZZ)$ with
\begin{equation}
\label{gammaaction}
\psi_{\gamma}(m,r,\ell)=(m+\ell\gamma, (d\alpha_{\gamma})(m)+\ell c_{\gamma}+r,
\ell).
\end{equation}
A basis of monomial 
sections of $\Gamma(X_{\Sigma}, \shL^{\otimes \ell})$ is indexed by
the set $C(\Delta_{\varphi})\cap (M\times\ZZ\times \{\ell\})$: for
$p$ in this set, we write $z^p$ for the corresponding section of
$\shL^{\otimes \ell}$. The above
action on $C(\Delta_{\varphi})$ lifts the $\Gamma$-action on $X_{\Sigma}$
to a $\Gamma$-action on each $\shL^{\otimes \ell}$. To write down sections
of $\shL^{\otimes \ell}$
on the quotient, one only need write down $\Gamma$-invariant sections
of $\shL^{\otimes \ell}$ on $X_{\Sigma}$, and
this can be done by taking, for any $m\in {1\over \ell}M$, the infinite sums
\[
\vartheta_m=\vartheta^{[\ell]}_m:=\sum_{\gamma\in\Gamma} z^{\psi_{\gamma}(\ell m,\ell\varphi(m),
\ell)}.
\]
We use the superscript $[\ell]$ to indicate the power of $\shL$ when
ambiguities can arise. We call $\ell$ the \emph{level} of the theta
function.

To see such an expression makes sense on the formal completion of the
zero fibre of $X_{\Sigma}\rightarrow\AA^1$, one focuses on an affine
chart of $X_{\Sigma}$ defined by a vertex $v=(m,\varphi(m))$ of 
$\Delta_{\varphi}$:
this affine chart is $\Spec \CC[T_v\Delta_{\varphi}\cap (M\oplus\ZZ)]$,
where $T_v\Delta_{\varphi}$ denotes the tangent cone to $\Delta_{\varphi}$
at the vertex $v$. One trivializes the line bundle $\shL^{\otimes \ell}$ 
in this chart using $z^{(\ell v,\ell)}\mapsto 1$, so that $\vartheta_m$ coincides
with the regular function
\[
\sum_{\gamma\in\Gamma} z^{\psi_{\gamma}(\ell m,\ell\varphi(m),\ell)
-(\ell v,\ell)}.
\]
Observe by the convexity of $\varphi$ that with $t=z^{(0,1)}
\in\CC[T_v\Delta_{\varphi}\cap (M\oplus\ZZ)]$, 
for any $k>0$ all but a finite number of monomials in this sum lie in
the ideal $(t^{k+1})$. Thus $\vartheta_m$ makes sense as a section
of $\shL^{\otimes \ell}$ on the $k$-th order thickening of the zero fibre
of $X_{\Sigma}\rightarrow \AA^1$, and since invariant under the $\Gamma$-action,
descends to a section on $\shA_k$.

Furthermore, one sees that $\vartheta_m=\vartheta_{m+\gamma}$, 
so if we set $B=M_{\RR}/\Gamma$, we obtain a set of theta functions
indexed by the points of $B({1\over \ell}\ZZ)$. One can show the
following facts:
\begin{enumerate}
\item The functions $\vartheta_m$ extend as holomorphic functions to
give the usual canonical theta functions on non-zero fibres of
$\pi^{-1}(D)/\Gamma\rightarrow D$.
\item The set $\{\vartheta_m \,|\, m\in B({1\over \ell}\ZZ)\}$ form a basis
for $\Gamma(\shA,\shL^{\otimes \ell})$ as a $\CC\lfor t\rfor$-module.
\item Denote by $\P$ also the polyhedral decomposition of
$B$ induced by the $\Gamma$-periodic decomposition of $M_{\RR}$. Assume no cells
of $\P$ are self-intersecting: this is equivalent to
all irreducible components of the central fibre
$\shA_0$ being normal. 
Each maximal cell $\sigma\in\P$ (thought of as a subset of
$B$) then corresponds to an irreducible component of the central fibre 
$\shA_0$ isomorphic to $\PP_{\sigma}$,
the projective toric variety determined by the lattice polytope $\sigma$.
Then if $m\in\sigma$, the restriction of $\vartheta_m$ to $\PP_{\sigma}$
is precisely the section of $\O_{\PP_{\sigma}}(\ell)\cong
\shL^{\otimes \ell}|_{\PP_{\sigma}}$ determined by $\ell m\in \ell\sigma$.
If $m\not\in\sigma$, then $\vartheta_m|_{\PP_{\sigma}}=0$. Thus
$\vartheta_m$ can be viewed as a lifting of the natural monomial section
of $\Gamma(\shA_0,\shL^{\otimes \ell}|_{\shA_0})$ which is non-zero on
those irreducible components indexed by $\sigma\in\P$ with $m\in\sigma$
and is given by $\ell m\in \ell\sigma$ on those components.
\end{enumerate}

This construction is particularly easy to describe as there is a global
description coming from the universal cover $M_{\RR}
\rightarrow B$. For more general $B$, we shall not have such a nice global
construction, and as a consequence, it is beneficial to give here a
more local description of theta functions. 

We can in fact use the $\Gamma$-action \eqref{gammaaction}
on $M\oplus\ZZ\oplus\ZZ$ to define a local system with
fibres $M\oplus\ZZ\oplus\ZZ$ on $B$ which we shall
call $\widetilde{\shP}$. The monodromy of the local system
is given by \eqref{gammaaction}; this uniquely determines the local system.
One checks one has the following commutative diagram of local systems on $B$:
\begin{equation}
\label{bigcommdiag}
\begin{split}
\xymatrix@C=30pt
{&0\ar[d]&0\ar[d]&&\\
&\underline{\ZZ}\ar[d]\ar[r]^=&\underline{\ZZ}\ar[d]&&\\
0\ar[r]&\shP\ar[d]^r\ar[r]&\widetilde\shP\ar[d]^{\tilde r}
\ar[r]^{\deg}&\underline{\ZZ}\ar[d]^{=}\ar[r]&0\\
0\ar[r]&\Lambda\ar[r]\ar[d]&\shAff(B,\ZZ)^*\ar[r]\ar[d]&\underline{\ZZ}\ar[r]
&0\\
&0&0&&
}
\end{split}
\end{equation}
Here $\underline{\ZZ}$ denotes the constant local system with stalks $\ZZ$,
the map $\deg$ is induced by the projection $M\oplus\ZZ\oplus\ZZ
\rightarrow\ZZ$ onto the last component, and $\shP$ is defined to be
the kernel of this map (hence has monodromy given by the restriction
of the action \eqref{gammaaction} to the first two components
$M\oplus\ZZ$). The inclusions
of $\underline{\ZZ}$ in $\shP$ and $\widetilde\shP$ are induced by
the inclusions $\ZZ\rightarrow M\oplus\ZZ, M\oplus\ZZ\oplus\ZZ$
into the second component. Here the quotient $\shP/\underline{\ZZ}$
is just the constant sheaf $\underline{M}$, and since $M_{\RR}$ is
canonically the tangent space to any point of $B$, we identify $\underline{M}$
with $\Lambda$, the local system of integral vector fields on $M$.
Finally, $\shAff(B,\ZZ)$ denotes the local system of integral affine linear
functions on $B$ (functions with integral slope and integral constant part),
and $\shAff(B,\ZZ)^*$ denotes the dual local system. To see the identification
of $\widetilde\shP/\underline{\ZZ}$ with $\shAff(B,\ZZ)^*$, write
$\Aff(M,\ZZ)=N\oplus\ZZ$, with $(n,c)\in N\oplus\ZZ$ defining the affine
linear map $m\mapsto \langle n,m\rangle+c$. Thus $M\oplus\ZZ$ is
canonically $\Aff(M,\ZZ)^*$. The action of $\Gamma$ on $\Aff(M,\ZZ)$
via pull-back of affine linear functions
is given by $(n,c)\mapsto (n,c+\langle n,\gamma\rangle)$, and the 
transpose action on $\Aff(M,\ZZ)^*$ is then precisely the restriction
of \eqref{gammaaction} to the first and third components of
$M\oplus\ZZ\oplus\ZZ$.

We can then describe a theta function as follows. Let $m\in B({1\over \ell}
\ZZ)$, and we want to describe $\vartheta_m$. We previously described 
$\vartheta_m$ as a sum of monomials $z^p$ with 
$p\in M\times \ZZ\times \{\ell\}$. Choose a point $x\in B$.
We can identify $\widetilde\shP_x$
with $M\oplus\ZZ\oplus\ZZ$ by choosing a lift $\tilde x\in M_{\RR}$
of $x$, so we can identify $\vartheta_m$ with a sum of monomials
$z^p$ with $p\in\widetilde\shP_x$. Note that the choice of lifting
is irrelevant, as a different lifting gives an identification related
by the transformation $\psi_{\gamma}$, and $\vartheta_m$ is invariant
under the action of $\psi_{\gamma}$.

We can then write
\begin{equation}
\label{pathdescriptiontheta}
\vartheta_m=\sum_{\delta} \Mono(\delta),
\end{equation}
where we sum over all affine linear maps $\delta:[0,1]\rightarrow B$
with the property that $\delta(0)=m$ and $\delta(1)=x$. We define
$\Mono(\delta)$ as follows. We have a canonical element of 
the stalk of $\shAff(B,\ZZ)^*$ at $m$ given by $\ell\cdot \ev_m$, where
$\ev_m$ denotes evaluation of integral affine functions at the point
$m$. This element can then be lifted to the stalk of $\widetilde\shP_m$
in a canonical way determined by $\varphi$: choosing a lifting $\tilde m$
of $m$ to $M_{\RR}$, we take $(\ell\tilde m,\ell\varphi(\tilde m),\ell)
\in M\oplus\ZZ\oplus\ZZ$; this defines a well-defined element of
$\widetilde\shP_m$ independent of the choice of lift $\tilde m$. This element
is indeed a lift of $\ell\cdot \ev_m$.
We call this element $m_{\varphi}\in \widetilde\shP_m$;
note that by construction $\deg(m_{\varphi})=\ell$. So far this is
independent of the choice of $\delta$. But now parallel transport 
$m_{\varphi}$ along the path $\delta$ to get an element 
$m_{\varphi}^{\delta}\in
\widetilde\shP_x$, and define $\Mono(\delta)=z^{m_{\varphi}^{\delta}}$.

It is not hard to check that \eqref{pathdescriptiontheta}
then coincides with our original description of $\vartheta_m$.

This does not represent anything radical: we are simply reinterpreting
the action of $\Gamma$ which led to theta functions in terms of the fundamental
group of $B$ in the guise of different choices of paths between
$m$ and $x$ (there is one such linear path for every choice of lift
$\tilde x$ of $x$ to $M_{\RR}$). However, now the description of theta
functions will generalize nicely. In particular, we can see the connection
between theta functions and homological mirror symmetry in a more direct manner.

To see this, we define a map
\[
{\bf vect}:\shAff(B,\ZZ)^*\rightarrow \shT_B
\]
as follows. An element of $\shAff(B,\ZZ)_x^*$ is an integral linear functional
on the vector space $\shAff(B,\RR)_x$ of germs of affine linear functions
at $x$ (with no integrality restriction). Restricting to the subspace
of functions which vanish at $x$, one obtains a derivation, yielding
a tangent vector at $x$. This defines the map.

For example, at $m\in B({1\over \ell}\ZZ)$, ${\bf vect}(\ell\cdot \ev_m)=0$
in $\shT_{B,m}$, simply because $\ev_m$ evaluates functions at $m$.
However, if $\delta:[0,1]\rightarrow B$ is a path with $\delta(0)=m$,
let $\tilde\delta:[0,1]\rightarrow M_{\RR}$ be a lifting with
$\tilde\delta(0)=\tilde m$. Let $\alpha(t)$ denote the
parallel transport of $\ell\cdot \ev_m$ along $\delta$ to
$\shAff(B,\ZZ)^*_{\delta(t)}$. Define
\begin{equation}
\label{displacementvectordef}
{\bf v}(t)={\bf vect}(\alpha(t)).
\end{equation}
Then one calculates that ${\bf v}(t)$ is the
tangent vector $\ell(\tilde m-\tilde\delta(t))$. In particular, provided 
$\delta$ is in fact linear, ${\bf vect}$ applied to the parallel transport
of $\ell\cdot\ev_m$ provides a vector field along $\delta$ which is
always tangent to the path $\delta$, always points towards the
initial point of the path, and increases in length as we move away
from the initial point at a rate proportional to $\ell$. 

The vector field ${\bf v}(t)$ gives rise to a holomorphic triangle
in $X(B)$ via
\begin{align*}
\psi:[0,1]\times[0,1] & \rightarrow X(B)\\
(t,s) &\mapsto s\cdot {\bf v}(t)\in \shT_{B,\delta(t)}\mod\Lambda_{\delta(t)}.
\end{align*}
Note this map contracts the edge of the square $\{0\}\times[0,1]$,
giving the triangle. This triangle is depicted in Figure 
\ref{sampletriangle}. 
Here $L_x$ is the fibre $\shT_{B,x}/\Lambda_x$ of the SYZ fibration
$X(B)\rightarrow B$.
This triangle can be seen as a contribution to the Floer multiplication
\[
\mu_2:\Hom(L_\ell,L_x)\times \Hom(L_0,L_\ell)\rightarrow\Hom(L_0,L_x).
\]
In particular, it yields a contribution to the product
of $p\in L_0\cap L_\ell$, corresponding to the point $m\in B({1\over \ell}\ZZ)$,
with the unique point of $L_\ell\cap L_x$.

This can be interpreted on the mirror side using homological mirror symmetry,
where we assume $L_\ell$ corresponds to $\shL^{\otimes \ell}$ and $L_x$ to
the structure sheaf of a point, as the composition map
\[
\Hom(\shL^{\otimes \ell},\O_x)\otimes
\Hom(\O_{\shA},\shL^{\otimes \ell})
\rightarrow\Hom(\O_{\shA},\O_x).
\]
Here $m$ determines the theta function $\vartheta_m\in\Hom(\O_{\shA},
\shL^{\otimes \ell})$ and a non-zero element of $\Hom(\shL^{\otimes \ell},\O_x)$
can be interpreted as specifying an identification $\shL^{\otimes \ell}\otimes
\O_x\cong \O_x$. The composition of $\vartheta_m$ with this
identification can then be viewed as specifying the value of the
section $\vartheta_m$ at the point $x$. Thus the description of
$\vartheta_m$ as a sum over paths $\delta$ then corresponds, naturally,
via this association of triangles to paths $\delta$, to the Floer theoretic
description on the mirror side.

\begin{figure}
\input{sampletriangle.pstex_t}
\caption{}
\label{sampletriangle}
\end{figure}

It is also important to describe multiplication of theta functions. 
Indeed, this allows us to describe the homogeneous coordinate ring
$\bigoplus_{\ell\ge 0} H^0(\shA,\shL^{\otimes \ell})$.
Given $m_i\in B({1\over \ell_i}\ZZ)$,
$i=1,2$, we wish to describe the coefficients of the expansion
\begin{equation}
\label{thetamultformula}
\vartheta_{m_1}\cdot \vartheta_{m_2} =\sum_{m\in B({1\over \ell_1+\ell_2}\ZZ)}
c_{m_1,m_2,m}\vartheta_m.
\end{equation}
It is not difficult to see that the coefficients are given by
\[
c_{m_1,m_2,m}=\sum_{\delta_1,\delta_2} t^{c(\delta_1,\delta_2)}
\]
where we sum over all straight lines $\delta_1,\delta_2:[0,1]\rightarrow
B$ connecting
$m_1, m_2$ to $m$ respectively, with the property that, if ${\bf v}_1$,
${\bf v}_2$ are defined by \eqref{displacementvectordef} using $\delta_1,
\delta_2$ respectively, then ${\bf v}_1(1)+{\bf v}_2(1)=0$.
We leave it to the reader to determine the exponent $c(\delta_1,\delta_2)
\in\NN$, depending on $\delta_1,\delta_2$: see \cite{DBr}, page 625,
in the case of the elliptic curve.

Again, this description of multiplication can be interpreted in terms
of Floer homology. Each pair $\delta_1,\delta_2$ contributing to the
sum gives rise to a triangle as depicted in Figure \ref{p1p2triangle}. 
Here, the triangle
is a union of two triangles in $X(B)$, fibering over $\delta_1$, $\delta_2$,
as depicted in that figure. The triangle on the left is determined
by $\delta_1$ as before, while the triangle on the right is the image
of the map
\begin{align*}
\psi:[0,1]\times [0,1]&\rightarrow X(B)\\
(t,s)&\mapsto \sigma_{\ell_1}(\delta_2(t))+s{\bf v}_2(t),
\end{align*}
where $\sigma_\ell$ is given by \eqref{sectiondef}.
The fact that these two triangles match up along the dotted line, which
lies over $m$, is just the statement that ${\bf v}_1(1)+{\bf v}_2(1)=0$.
The number $c(\delta_1,\delta_2)$ can be seen to be related (but not
equal to) the symplectic area of this triangle, again see \cite{DBr},
pp.\ 626--628. 

The multiplication formula \eqref{thetamultformula} 
can be viewed as a kind of global generalization
of a much simpler rule for multiplying sections of powers of a given
ample line bundle on a toric variety. Indeed, such a line bundle $\shL$ 
determines
a lattice polytope $B \subset \RR^n$, and the points of $B({1\over \ell}\ZZ)$
correspond to a monomial basis for the global sections of $\shL^{\otimes \ell}$.
The product of the sections corresponding to $m_i\in B({1\over \ell_i}\ZZ)$,
$i=1,2$, is just the section corresponding to the weighted average
$m=(\ell_1m_1+\ell_2m_2)/(\ell_1+\ell_2)\in B({1\over \ell_1+\ell_2}\ZZ)$. 
This can be interpreted
in terms of paths $\delta_1$, $\delta_2$ joining $m_1$ and $m_2$ to $m$,
as in \eqref{thetamultformula}.

\begin{figure}
\input{p1p2triangle.pstex_t}
\caption{}
\label{p1p2triangle}
\end{figure}

\section{Singularities, theta functions, jagged paths}

We now would like to generalize these constructions to affine manifolds
with singularities, as is necessary if we are to obtain any interesting
examples. We will begin with some simple examples to provide guidance.
In particular, we will take $B$ to be a compact affine manifold with
boundary, analogous to the very simple case where $B$ is just a lattice
polytope, but allow a few simple singularities to appear in $B$.
We always assume $B$ is locally convex along $\partial B$.

\subsection{The basic example}
\label{basicexamplesection}
\begin{figure}
\input{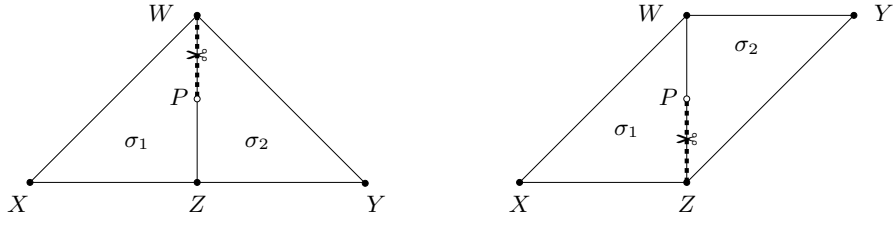}
\caption{The affine manifold $B_1$. The diagram shows the affine embeddings
of two charts, obtained by cutting the union of two triangles as indicated
in two different ways. Each triangle is a standard simplex.}
\label{TwoTriangles}
\end{figure}

We will revisit some examples introduced in \cite{GSInv}.
For now, we consider the simplest example,
the affine manifold $B_1$ given in Figure \ref{TwoTriangles}.
The points of $B_1(\ZZ)$ are labelled in the diagram as $X,Y,Z$ and $W$.
The affine structure has one singularity, the point $P$, and the point
$P$ can be chosen freely within the line segment joining $W$ and $Z$.
The piecewise linear function $\varphi$ takes the value $0$ at $X,W$ and $Z$
and $1$ at $Y$.

If we were in the purely toric case, say $B_1$ being either the polygon
pictured on the left or the right in Figure \ref{TwoTriangles}, then
each integral point would represent a purely monomial section of the line
bundle on the toric variety corresponding to $B_1$. Further,
the multiplication law for monomials would be either $XY=Z^2$ (on the
left) or $XY=WZ$ (on the right). Each such product can be viewed
by giving paths $\delta_1, \delta_2$ with initial endpoints $X$ and $Y$
respectively and terminating at either the point $Z$ or the point
${1\over 2 }(W+Z)$ (viewing these points as elements of $\RR^2$ rather than
as variables). If we also took into account the polyhedral decomposition
and the choice of $\varphi$, one obtains a degeneration of one of these
two toric varieties into a union of two planes, given by the equation
$XY=tZ^2$ or $XY=tWZ$ in the two cases, where $t$ is the deformation parameter.

However, we are not in the purely toric case, and if we follow
the philosophy of the previous section, the description of the sections
specified by the points of $B_1(\ZZ)$ and their multiplication rule should
be determined by drawing straight lines. Let us first consider 
heuristically how lines should contribute to the description of sections, 
and then consider how we should think of the product rule.

First looking at the points labelled $W$ or $Z$, we note that given any
reference point $x\in B_1\setminus\{P\}$, 
there is a unique line segment joining $W$ or $Z$ to 
$x$. In analogy with the abelian variety case, we would expect this to
tell us that the corresponding sections are still represented by monomials
at $x$. Next, consider the point $X$. If $x$ is contained
in $\sigma_1$, then there is again a unique line segment joining $X$ and
$x$, so we expect a monomial representative for this section. On the
other hand, if $x$ lies in $\sigma_2$, we may have either one or two
straight lines, depending on the precise location of $x$: In Figure
\ref{TwoTriangles3}, there is a line joining $X$ and $x$ in the first
chart but not in the second, as the line drawn in the right-hand chart crosses
the cut.
On the other hand, in Figure \ref{TwoTriangles4}, there are in fact
two distinct line segments joining $X$ and $x$.

One solution to this ambiguity is to simply include in the sum a contribution
from the line segment in the right-hand chart of Figure \ref{TwoTriangles3}.
However, this is not a straight line: if drawn in the correct chart,
it becomes bent, as depicted in Figure \ref{TwoTriangles5}.

How do we justify counting this bent line, or
as we shall call it, \emph{jagged path}? The explanation is that when we 
introduce singularities, as explained in \cite{GSInv}, we need
to introduce walls emanating from the singularity, in this case rays
heading in the direction of the points $W$ and $Z$. Each ray has
a function attached to it; the precise role that this function plays
will be explained later. But the essential point is that we no longer
need to use straight lines to join $X$ and $x$. We will allow our lines
to bend in specified ways when lines cross walls. In the case of $B_1$,
this bending exactly accounts for the jagged path in Figure \ref{TwoTriangles5}.
As a consequence, we should expect the section
corresponding to $X$ will be represented as a sum of two monomials
in any event in a chart corresponding to the right-hand side of
$B_1$, regardless of the position of $x$.

As we have not yet been very clear what these charts mean and how we
are representing sections in general, it is perhaps more informative to
loook at the product $XY$. 
This product should be given as a sum over all suitable
choices of paths $\delta_1, \delta_2$. To realise this,
let us first assume the point $P$ lies below ${1\over 2}(W+Z)$. Then
Figure \ref{TwoTriangles2} shows two possible choices of the pairs
$\delta_1,\delta_2$, and the product $XY$ should be determined as a sum
over these two ways of averaging $X$ and $Y$: the presence of the
singularity has created this ambiguity. 
Taking the PL function $\varphi$ into account, the
suggestion then is that we should have the multiplication rule
\[
XY=t(Z^2+WZ).
\]
Note that this gives a family over $\Spec \CC[t]$ whose fibre over $t=0$
is a union of two $\PP^2$'s in $\PP^3$
(determined by the two standard simplices
$\sigma_1$ and $\sigma_2$), and for general $t$, we obtain a non-singular
quadric surface in $\PP^3$.

\begin{figure}
\input{TwoTriangles3.pstex_t}
\caption{}
\label{TwoTriangles3}
\end{figure}

\begin{figure}
\input{TwoTriangles4.pstex_t}
\caption{}
\label{TwoTriangles4}
\end{figure}

\begin{figure}
\input{TwoTriangles5.pstex_t}
\caption{}
\label{TwoTriangles5}
\end{figure}

\begin{figure}
\input{TwoTriangles2.pstex_t}
\caption{}
\label{TwoTriangles2}
\end{figure}

Now this argument depended on the fact that the point $P$ was chosen
below the half-integral point ${1\over 2}(W+Z)$.
Since there is no sense that this singular
point has a natural location, this is not particularly satisfactory. If we
move $P$ above this point, the second choice of $\delta_1,\delta_2$
seems to disappear. 

The solution again is to use the walls to provide corrections. If
$P$ lies above the point ${1\over 2}(W+Z)$, then we will again obtain two
contributions as depicted in Figure \ref{TwoTriangles6}, where this time
a piece of the wall emanating from $P$ is used to correct for the fact
that for the labelled $\delta_1$, $\delta_2$, we do not have ${\bf v}_1(1)
+{\bf v}_2(1)=0$ (where ${\bf v}_i$ is defined using 
\eqref{displacementvectordef}). Rather, we have ${\bf v}_1(1)+{\bf v}_2(1)+w
=0$, where $w$ is the unit tangent vector pointing from $P$ to $W$.

These pictures can be justified heuristically in terms of holomorphic disks
contributing to Floer multiplication. 
If we think of a space $X_1$ fibering in tori
over $B_1$, there is a singular fibre over $P$, a two-torus with a circle
pinched to a point. We expect that $X_1$ will contain holomorphic disks
fibering over the two walls emanating from $P$. In particular, for any
point $y$ in $B_1\setminus\{P\}$ on the line 
segment $\overline{WZ}$, there is a holomorphic
disk in $X_1$ with boundary contained in the fibre over $y$. We can use
this to build a piecewise linear disk as follows.

Let $\delta:[0,1]\rightarrow B_1$ be the parameterized jagged path of
Figure \ref{TwoTriangles5}, bending at time $t_0\in (0,1)$. 
Modify the definition of ${\bf v}:[0,1]\rightarrow \shT_{B_1}$ defined
using \eqref{displacementvectordef} by taking it to coincide with
${\bf v}$ of \eqref{displacementvectordef} for $0\le t<t_0$ and with
${\bf v}+w$ for $t_0\le t\le 1$. One checks that ${\bf v}(t)$ is always
tangent to $\delta$.

As in \S\ref{abeliansection}, we use ${\bf v}$ to define a polygon in
$\shT_{B_1}$, but because of the discontinuity in ${\bf v}$ we obtain
a picture as in Figure \ref{PLdisk}. Recall that $X(B_1\setminus\{P\})$
is obtained by dividing the tangent spaces of $B_1\setminus \{P\}$ 
by integral vector fields. If ${\bf w}_{\pm}
=\lim_{t\rightarrow t_0^\pm}{\bf v}(t)$, 
then ${\bf w}_+-{\bf w}_-=w$, so in fact
${\bf w}_+$ and ${\bf w}_-$ are identified
in $X(B_1\setminus\{P\})$. So the line segment joining ${\bf w}_+$ and
${\bf w}_-$ becomes a loop. We can then glue in a holomorphic disk
emanating from $P$, attaching its boundary to the loop. This gives
the triangle which is, roughly speaking, the contribution to Floer 
multiplication describing the section $X$.

A similar picture explains the contributions to the product $XY$:
the failure of ${\bf v}_1(1)+{\bf v}_2(1)=0$ is dealt with by gluing
in the holomorphic disk.

This is just a heuristic: these disks are not actual holomorphic disks.
However, a variant of this example is considered in great detail in
\cite{P11} and the result here agrees with the actual result from
Floer multiplication. The advantage for us is that we can describe
everything combinatorially.

Before looking at some more complex examples, let us give a more precise
description of what we are doing. Unfortunately, doing so requires a number
of technical details; we shall try to avoid the most unpleasant aspects,
but the reader should be advised in what follows that the definitions
are only approximately correct!

\begin{figure}
\input{TwoTriangles6.pstex_t}
\caption{}
\label{TwoTriangles6}
\end{figure}

\begin{figure}
\input{PLdisk.pstex_t}
\caption{}
\label{PLdisk}
\end{figure}

\subsection{Structures and jagged paths}

Before we get into details, suppose we are given data $(B,\P,\varphi)$,
where $B$ is an integral affine manifold with singularities, $\P$ is a 
polyhedral
decomposition of $B$, and $\varphi$ is an integral multi-valued PL function
on $B$. We then obtain a generalization of the diagram \eqref{bigcommdiag}
of sheaves on $B_0$. This is a direct generalization of the torus case.
Saying that $\varphi$ is a multi-valued PL function on $B$ means
that there is an open cover $\{U_i\}$ of
$B_0$ such that $\varphi$ is represented by a single-valued function
$\varphi_i$ on $U_i$, with $\varphi_i-\varphi_j$ affine linear on
$U_i\cap U_j$. We then construct the local system $\widetilde\shP$
on $B_0$ as follows.
On $U_i$, $\widetilde{\shP}$ is isomorphic to
$\underline{\ZZ}\oplus \shAff(B_0,\ZZ)^*|_{U_i}$, and on
$U_i\cap U_j$, $\underline{\ZZ}\oplus\shAff(B_0,\ZZ)^*|_{U_i}$
is identified with
$\underline{\ZZ}\oplus\shAff(B_0,\ZZ)^*|_{U_j}$ via
\[
(\ell,\alpha)\mapsto (\ell+\alpha(\varphi_j-\varphi_i),\alpha)
\]
for $\alpha\in\Gamma(U_i\cap U_j,\shAff(B_0,\ZZ)^*)$ and $\ell\in \ZZ$.
We then have the projection map $\tilde r:\widetilde{\shP}\rightarrow
\shAff(B_0,\ZZ)^*$, and dualizing the exact sequence
\[
0\rightarrow\underline{\ZZ}\rightarrow\shAff(B_0,\ZZ)\rightarrow
\check\Lambda\rightarrow 0
\]
(with the third arrow given by exterior derivative) gives the bottom
row of \eqref{bigcommdiag}. From this follows the whole diagram,
with $\shP$ defined as the kernel of the map $\deg$.

In \cite{GS11}, we gave the definition of a \emph{structure} for
an integral affine manifold with singularities $(B,\P)$. 
Structures were used for the proof of Theorem 
\ref{mainsmoothingtheorem} to encode the explicit data necessary
to describe the smoothing. 
For $B$ of 
arbitrary dimension, a structure $\foD$ is a collection of \emph{slabs} and
\emph{walls}. These are codimension 1 polyhedra in $B$ which are either
contained in codimension one cells of $\P$ (slabs) or contained in
maximal cells of $\P$ but not contained in codimension one cells (walls). 
Slabs and walls carry additional data, certain formal power series which
are used to describe gluing automorphisms. In order for the gluing to
be well-defined, a structure must satisfy the notion of \emph{compatibility}.
Producing a compatible structure is the main work of \cite{GS11}. This procedure
is described at greater length in \cite{GSInv} and is covered in full
detail in \cite{Gr11}, Chapter 6 in the two-dimensional case. Slabs and walls
need to be treated somewhat differently in the algorithm of \cite{GS11}
for producing compatible structures for technical reasons, but in the
context here we will essentially be able to ignore these issues.

A number of details of this construction
are surveyed in \cite{GSInv}. The crucial points to know are the following:
\begin{enumerate}
\item The deformation $\check\shX\rightarrow\Spec \CC\lfor t\rfor$ is 
given order-by-order, with $\check\shX_k:=\check\shX\times_{\CC\lfor t\rfor}
\CC[t]/(t^{k+1})$ constructed explicitly from a structure.
\item $\check\shX_k$ is a thickening of $\check\shX_0$, and hence has the same set
of irreducible components, indexed by $\sigma\in\P_{\max}$. The function
$\varphi$ determines some standard toric thickenings of affine open subsets
of $\PP_{\sigma}$ for $\sigma\in\P_{\max}$. Then $\check\shX_k$ is obtained
by using the structure to glue together these standard pieces in non-toric
ways.
\item If $x\in B_0$ there is
a monoid $P_x\subseteq \shP_x$, defined using $\varphi$,
along with an inclusion $\NN\rightarrow P_x$, yielding a family
$\pi_x:\Spec \CC[P_x]\rightarrow \Spec \CC[t]$. This provides a local
model for the smoothing, as follows.
Let $\tau\in\P$ be the smallest
cell containing $x$. There is a one-to-one correspondence between
cells $\sigma\in\P_{\max}$ containing $x$ and irreducible components of
$\pi_x^{-1}(0)$. Furthermore, the irreducible component corresponding to
$\sigma$ is isomorphic to the affine open subset of $\PP_{\sigma}$
determined by the face $\tau$ of $\sigma$.
Then the corresponding irreducible component of $\Spec\CC[P_x]\times
_{\AA^1} \Spec \CC[t]/(t^{k+1})$ is the standard toric thickening of
this affine open subset of $\PP_{\sigma}$.
\item In \cite{GHKS}, we will explain how the line bundle $\shL^{\otimes d}$ 
is obtained by gluing standard
line bundles on these standard pieces, again with gluing dictated by
the structure. For any given $\ell$, there is a $P_x$-torsor $Q_x^\ell
\subseteq \widetilde\shP_x\cap \deg^{-1}(\ell)$ defining a line bundle
on $\Spec\CC[P_x]$; its restriction to the irreducible components of
$\Spec\CC[P_x]\times _{\AA^1} \Spec \CC[t]/(t^{k+1})$ describes these
standard line bundles.
\end{enumerate}
The definition of the monoid $P_x$ is discussed in \cite{GSInv}, \S3.1; 
the details will not be so important here. The $P_x$-torsor $Q^\ell_x$ has
not yet been discussed in the literature.

As most of the conceptual issues are already present in the two-dimensional
case, instead of using walls and slabs, we can use rays, much as 
\cite{KS06} had done. For precise details of what follows, see \cite{Gr11},
Chapter 6.

Assume $B$ is two-dimensional. Roughly, a ray consists of:
\begin{enumerate}
\item a parameterized ray or line segment $\fod:[0,\infty)\rightarrow B$
or $\fod:[0,1]\rightarrow B$, with image a straight line of rational slope.
A ray will continue until it hits a singularity or the boundary of $B$;
otherwise it continues indefinitely.
\item A formal power series
\[
f_{\fod}=1+\sum_p c_pz^p
\]
where $p$ runs over the set of global sections of $\fod^{-1}\shP$
with the property that $r(p)$ is always tangent to the image of
$\fod$, \emph{negatively} proportional to the derivative $\fod'$. There are
more constraints on $f_{\fod}$ necessary to guarantee convergence of 
the gluing construction considered below, but we won't worry about these
technical details as all examples considered here will be quite simple.
\end{enumerate}

A \emph{structure} $\foD$ is then a collection of rays $\{(\fod,f_{\fod})\}$.

\begin{example}
\label{basicexamplescatteringdiagram}
In the basic example of \S\ref{basicexamplesection}, the relevant
structure consists of two rays: $\fod_1$, a line segment from $P$ to $W$,
and $\fod_2$, a line segment from $P$ to $Z$. Since $\varphi$ is single-valued,
we can write $\shP=\underline{\ZZ}\oplus\Lambda$, and we take
\[
f_{\fod_1}=1+z,\quad f_{\fod_2}=1+w,
\]
where $z$ and $w$ are the monomials corresponding to 
$(0,(0,-1))$ and $(0,(0,1))$ respectively. Note applying $r$
to these two elements gives the primitive tangent vectors to the segment
$\overline{ZW}$ pointing towards $Z$ and $W$ respectively. 
We shall by abuse of
notation refer to these tangent vectors as $z$ and $w$ also.
\qed
\end{example}

Note that a point $m\in B({1\over \ell}\ZZ)$ defines an element $\ell\cdot
\ev_m$ of $\shAff(B_0,\ZZ)^*_m$ as in \S\ref{abeliansection}. In addition,
choosing
a local representative $\varphi_m$ for $\varphi$ in a neighbourhood
of $m$ gives a splitting $\widetilde{\shP}=\underline{\ZZ}\oplus
\shAff(B_0,\ZZ)^*$ in a neighbourhood of $m$, and we define an element
$m_{\varphi}\in\widetilde{\shP}_m$ by 
\[
m_{\varphi}=(\ell\cdot\varphi_m(m), \ell\cdot \ev_m)\in
\ZZ\oplus\shAff(B_0,\ZZ)^*_m.
\]
One notes this is independent of the choice of $\varphi_m$.
In fact, $m_{\varphi}$ lies in the $P_m$-torsor $Q^\ell_m$, and hence
defines a local section of the line bundle $\shL^{\otimes \ell}$. 

For example, in the case of $B=B_1$, the points $X,Y,Z$ and $W$
define elements of the stalks $\widetilde{\shP}_X,\ldots,\widetilde{\shP}_W$
of degree $1$. However, as parallel transport of $X$ and $Y$
is not well-defined because of monodromy of the local system
$\widetilde\shP$ around $P$, 
we cannot directly view these as defining global sections of $\shL$.

To do so requires the precise notion of jagged path. This definition
will appear in \cite{GHKS}.

\begin{definition}
\label{jaggeddef}
A \emph{jagged path in $B$ with respect to a structure $\foD$} consists of
the following data:
\begin{itemize}
\item[(a)] A continuous piecewise linear path $\gamma:[0,1]\rightarrow B$.
\item[(b)] For every maximal domain of linearity $L\subseteq [0,1]$ of
$\gamma$ we are given a monomial 
\[
m_L=c_Lz^{q_L}\in \CC[\Gamma(L,\gamma^{-1}(\widetilde{\shP})|_L)]
\]
\end{itemize}
satisfying the following two properties:
\begin{enumerate}
\item If $t\in (0,1)$ is a point contained in the interior of $L$
a maximal domain of linearity, then $\gamma'(t)$ is negatively proportional
to ${\bf vect}(\tilde r(q_L))$, where $\tilde r$ is defined in 
\eqref{bigcommdiag}. 
\item Let $t\in (0,1)$ be a point at which $\gamma$ is not affine linear,
passing from a domain of linearity $L$ to $L'$, with $y=\gamma(t)$.
Let $\{(\fod_j,x_j)\}$ be the set of pairs $\fod_j\in\foD$, $x_j$ in
the domain of $\fod_j$ such that $\fod_j(x_j)=y$. Let $n_{\fod_j}
\in \shAff(B_0,\ZZ)_y$ be the germ of a primitive integral affine linear
function which vanishes on the image of $\fod_j$ near $y$ and is positive
on the image of $L$ near $y$. We assume that $n=n_{\fod_j}$ can be chosen
independently of $j$; this is an assumption on the genericity of
$\gamma$ and can always be achieved by perturbing the endpoint $\gamma(1)$.
Expand
\begin{equation}
\label{bendingexpansion}
\prod_j f_{\fod_j}^{\langle n,\tilde r(q_L)\rangle}
\end{equation}
as a sum of monomials with distinct exponents. Note each monomial can be
viewed as an element of $\CC[\shP_y]\subseteq\CC[\widetilde\shP_y]$.

Then there is a term
$cz^q$ in this expansion with
\[
m_{L'}=m_L\cdot (cz^q)=c_L c z^{q_L+q}.
\] 
\end{enumerate}
\end{definition}

There are several important features of this definition. First, item (1)
says the monomial attached to the line segment always tells us the direction
of travel of the line segment. Since (2) tells us how these monomials change
at bends, it also tells us precisely how jagged paths bend.

Second, the exponents in \eqref{bendingexpansion} are by construction
always non-negative. This is important as $f_{\fod_j}$ need
not be invertible in the relevant rings. In fact, more naive approaches
to writing down sections of $\shL$ founder on precisely this point.

\begin{definition}
For $m\in B({1\over \ell}\ZZ)$ and $x\in B_0$ general, a \emph{jagged path} from
$m$ to $x$ is a jagged path satisfying
\begin{enumerate}
\item $\gamma(0)=m$;
\item $\gamma(1)=x$;
\item If $L$ is the first domain of linearity of $\gamma$, then
$m_L=z^{m_{\varphi}}$.
\end{enumerate}
\end{definition}

\begin{example} 
Let's examine this in detail with $B=B_1$, $\foD$ as in Example 
\ref{basicexamplescatteringdiagram}. Take $m=X\in B_1(\ZZ)$.
If $x\in\Int(\sigma_1)$, there is one jagged path from $X$ to $x$,
which just serves to parallel transport $m_{\varphi}$ to $x$.

On the other hand, to be explicit, let's take $x$ to be the point with
coordinates $(1/8,1/4)$ in $\sigma_2$, in the left-hand chart of
Figure \ref{TwoTriangles}, assuming $Z$, $Y$ and $W$ have coordinates
$(0,0)$, $(1,0)$ and $(0,1)$ respectively. Suppose a jagged path from $X$ to
$x$ crosses the segment $\overline{WZ}$ below $P$ (which we will take to lie
at $(0,1/2)$ to be explicit). Then we are crossing the ray $\fod_2$,
and we take $n_{\fod_2}$ to be the linear function $(a,b)\mapsto -a$.
Then $\tilde r(m_{\varphi})=\ev_X$ takes the value $1$ on $n_{\fod_2}$,
precisely because $n_{\fod_2}$ takes the value $1$ at $X$. Thus 
\eqref{bendingexpansion} is just given by $1+w$. 
So we can take
$cz^q=1$, in which case there is no change to the monomial and no bend in
the jagged path; this just yields the parallel transport of $X$ into
$\sigma_2$. Otherwise, we take $cz^q=w$. This replaces $X$ with
$Xw$, and changes the direction by noting 
that 
${\bf vect}(\ev_X+w)={\bf vect}(\ev_X)+w$. 
Recall here that we are using
the same notation $w$ for the tangent vector and corresponding
monomial. In particular, at the point of
intersection of the segment $\overline{WZ}$ 
with the jagged path, say at $(0,h)$,
we have
\[
{\bf vect}(\ev_X)=(-1,0)-(0,h)=(-1,-h),
\]
so the new direction is $-(-1,-h+1)$. In order for this path to then
pass through $x$, we need to take $h=1/3$.

A bit of experimentation shows that as we move the point $x$ around
inside $\sigma_2$, we always have \emph{precisely two} jagged paths from
$X$ to $x$, and the sum of the final attached monomials is independent
of the location of $x$. However, it is possible that both such jagged
paths are in fact straight, when viewed in the correct charts, 
as happens in Figure \ref{TwoTriangles4}.
To describe this sum, let us continue to denote
by $X$ the monomial parallel transported
from $\sigma_1$ into $\sigma_2$ below $P$. (Parallel transport
is carried out in the local system $\widetilde\shP$). Then the sum is $X+Xw$. 
\end{example}

It is this invariance which is the crucial property of jagged paths. 
To get our hands on this invariance, we use:

\begin{definition}
\label{Def: Lift(m)}
Let $m\in B({1\over \ell}\ZZ)$. For a jagged path $\gamma$ from
$m$ to $x$, let $\Mono(\gamma)\in \CC[\widetilde{\shP}_x]$
be the monomial attached to the final domain of linearity of $\gamma$.
Let 
\[
\Lift_x(m):=\sum_{\gamma} \Mono(\gamma)
\]
be the sum over all distinct jagged paths from $m$ to $x$. Here we view two jagged
paths to be the same if they just differ by a reparametrization of their
domains.
\end{definition}

Given a path $\gamma$ in $B_0$ connecting two points $x_1$ and $x_2$,
we can define a transformation $\theta_{\gamma,\foD}$ which, roughly
speaking, is a map $\CC[\widetilde{\shP}_{x_1}]\rightarrow
\CC[\widetilde{\shP}_{x_2}]$ which is given by parallel transport
and a composition of wall-crossing automorphisms: if we 
cross a ray $\fod$ at time $t$, we apply the transformation
\[
z^q\mapsto z^q f_{\fod}^{\langle n_{\fod},\tilde r(q)\rangle}
\]
for $q\in\widetilde{\shP}_{\gamma(t)}$, $n_{\fod}$ a primitive
integral affine
linear function vanishing along the image of $\fod$ near $\gamma(t)$
and positive on $\gamma(t-\epsilon)$. 

\begin{definition}
A structure $\foD$ is \emph{consistent} if for any path $\gamma$ from
$x_1$ to $x_2$ general points in $B_0$, $m\in B({1\over \ell}\ZZ)$,
\[
\Lift_{x_2}(m)=\theta_{\gamma,\foD}(\Lift_{x_1}(m)).
\]
In particular if $\gamma$ does not cross any ray of $\foD$, then 
$\Lift_x(m)$ is invariant under parallel transport.
\end{definition}

If we have a consistent structure $\foD$, we can then use the lifts
$\Lift_x(m)$ to define a global section $\vartheta_m$ of $\shL^{\otimes \ell}$. 
We can write down well-defined local descriptions of the section on the
various affine pieces of the irreducible components of $\check\shX_k$, and
consistency then guarantees that these local descriptions glue.

A main result of \cite{GHKS} is then:

\begin{theorem}
\label{GHKSmain}
\begin{enumerate}
\item The compatible structures constructed in \cite{GS11} are in fact
consistent.
\item Given a compatible and consistent structure, giving a formal degeneration
$\foX\rightarrow\Spf\CC\lfor t\rfor$, the above construction gives for
every $m\in B({1\over \ell}\ZZ)$ a section $\vartheta_m^{[\ell]}$ 
of the line bundle
$\shL^{\otimes \ell}$. This section has the property that for any 
$\sigma\in\P_{\max}$, $\vartheta_m^{[\ell]}|_{\PP_{\sigma}}$ is $0$ if $m
\not\in\sigma$ and otherwise coincides with the monomial section of
$\O_{\PP_{\sigma}}(\ell)$ defined by $m$.
\end{enumerate}
\end{theorem}

These are our theta functions. We call $\vartheta_m^{[\ell]}$ a theta function
of \emph{level $\ell$}. Keeping in mind that $m$ can lie in 
$B({1\over \ell}\ZZ)$ for various $\ell$, $\vartheta_m^{[\ell]}$ may depend on 
the level. However, we will write $\vartheta_m$ when not ambiguous.

The proof of this result is a fairly straightforward extension of arguments
given in \cite{GHKI} and \cite{CPS}.

We can also use jagged paths to describe multiplication. In general, 
for $m_1\in B({1\over \ell_1}\ZZ)$, $m_2\in B({1\over \ell_2}\ZZ)$, 
we should have a multiplication rule
\begin{equation}
\label{thetamultformula2}
\vartheta_{m_1}^{[\ell_1]}\cdot\vartheta^{[\ell_2]}_{m_2}
=\sum_{m\in B({1\over \ell_1+\ell_2}\ZZ)} c_{m_1,m_2,m}\vartheta_m^{[\ell_1
+\ell_2]}.
\end{equation}
As before, the coefficient should be determined as a sum over
pairs of jagged paths $\delta_1,\delta_2$, with $\delta_i$ a jagged
path from $m_i$ to $m$, and with balancing ${\bf v}_1(1)+{\bf v}_2(1)=0$.
There is a slight subtlety which we have already seen in 
\S\ref{basicexamplesection}: since $m$ is not free to be chosen generally,
it may lie on a ray of $\foD$. Indeed, this happens in 
\S\ref{basicexamplesection}.
So it is possible that balancing fails, and this is corrected by using
a contribution from the ray. 
Essentially, this can be accomplished by perturbing the point $m$ a little
bit, so that one of $\delta_1$ or $\delta_2$ has a chance to have an
addditional bend along that ray. This needs to be done with a bit of care,
so we omit the details of this. 

\subsection{Additional examples}

In \cite{GSInv}, we considered a number of other two- and three-dimensional
examples. Here, we will describe their homogeneous coordinate rings in terms
of jagged paths.

\begin{example}
The affine manifold $B_2$ is depicted in two left-hand diagrams
in Figure \ref{B2figure}.
The required structure $\foD$ is exactly as in the case of $B_1$,
with two rays, one from $P$ to $W$ and one from $P$ to $U$. We write
$X$ instead of $\vartheta_X$ etc.\ for the theta functions of level $1$.
One sees that the products of theta functions $WY$ and $UZ$ each
correspond to the theta function of level $2$ 
determined by the barycenter of the
square $\sigma_2$, so one obtains the purely toric relation $WY-ZU=0$.
Much as in the case of $B_1$, one obtains products $XY=t(U^2+UW)$ and
$XZ=t(W^2+WU)$; the two choices of pairs of jagged paths contributing
to the latter product are indicated in the right-hand diagram of
Figure \ref{B2figure}, if we assume
that the singular point occurs below the midpoint of the segment
$\overline{WU}$.
\end{example}

\begin{figure}
\input{B2figure.pstex_t}
\caption{}
\label{B2figure}
\end{figure}

\begin{figure}
\input{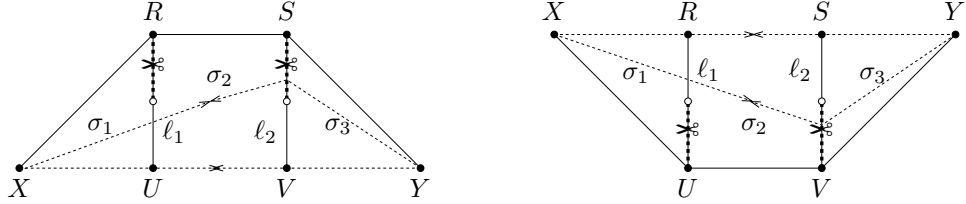}
\caption{Two of the pairs of jagged paths contributing to $XY$ are drawn
in each chart. Contrary to appearances, none of the eight jagged paths
appearing here bend!}
\label{B3figure2}
\end{figure}

\begin{example}
In Figure \ref{B3figure}, we have an example with two singularities,
with both charts shown. 
Note here we have four rays in the relevant structure,
two each emanating from the singular points. We 
take the functions attached to
the rays to be $1+u$ and $1+r$ as appropriate, where
$r$ and $u$ are the monomials corresponding to the tangent vectors
$(0,1)$ and $(0,-1)$ respectively.
We also take $\varphi$
to take the value $0$ on the square and the value $1$ at $X$ and $Y$.
One sees easily, using the same strategies as above,
various quadratic relations on the level $1$ theta functions. 
First, we have the purely
toric relation $RV=SU$, as neither of these products involve jagged paths
which cross rays. Next, the products $XS$, $XV$, $RY$ and $UY$ all behave
as in the previous example, and we can write
\[
XS=t(R^2+UR), \quad XV=t(U^2+UR), \quad RY=t(S^2+SV), \quad UY=t(V^2+VS).
\]
Finally, the product $XY$ is the most interesting, with four contributions,
\[
XY=t^2(UV+US+RV+RS).
\]
Figure \ref{B3figure2} shows all four pairs of jagged paths contributing
to these terms. Note in \cite{GSInv}, Example 3.4, this relation was obtained
from the previous ones by saturation of ideals.
\end{example}

\begin{figure}
\input{B3figure.pstex_t}
\caption{}
\label{B3figure}
\end{figure}

\begin{example}
We consider next $B_3$ of \S4.2 of \cite{GSInv}, see Figure \ref{B4figure}.
There are two singularities; the second diagram shows the affine
embedding of an open set containing the two cuts of the first chart.
Here $\varphi$ takes the value $0$ at $U, Z$ and $X$, and the value $1$
at $Y$ and $W$.

The structure $\foD$ defining the deformation is as follows. Each 
singularity produces two rays contained in the line segments containing
them; thus two of these rays intersect at $U$. Because of this intersection,
an extra ray, $\fop$ as drawn, is necessary to produce compatibility. The
function attached to $\fop$ is $f_{\fop}=1+t^2x^{-1}z^{-1}$, where 
$x$ and $z$ are the monomials corresponding to the tangent vectors $(1,0)$ and
$(0,1)$ respectively. The relations are then given by
\[
XY=t(U^2+UW), \quad ZW=t(U^2+YU).
\]
The term $U^2$ in the expression for $XY$ and $ZW$ is the expected
toric one from the first chart; the left-hand diagram of
Figure \ref{B4figure2} exhibits the
diagrams for the other two relations. In this figure we take the singularities
to be closer to the boundary of $B_3$ than to $U$.

More interesting in this case is how consistency arises. Let us consider
jagged paths from $Y$ into the chamber between the ray $\fop$ and
the segment $\overline{UX}$. In particular, consider a jagged path
that bends at $\fop$, e.g., $\delta_1$ as depicted in the right-hand
diagram of Figure \ref{B4figure2}.
If we take $U$ to be the origin, and if $\delta_1$ hits
$\fop$ at the point $(r,r)$, then the direction of the jagged path
after the bend is $(r+2,r+1)$, and hence always has slope $\ge 1/2$. Thus
such a jagged path can never end at a point of the chamber in question 
below the
line of slope $1/2$ passing through the origin. In order for the lift of
$Y$ to be independent of the endpoint inside this chamber, there must be some 
other jagged paths contributing the same monomial at points below the slope
$1/2$ line. One sees that a path of type $\delta_2$ does the trick. Note
also that if we had omitted the ray $\fop$, the jagged path $\delta_2$
would still exist, but consistency would fail because there would be no
substitute for $\delta_1$ above the slope $1/2$ line.
\end{example}

\begin{figure}
\input{B4figure.pstex_t}
\caption{}
\label{B4figure}
\end{figure}

\begin{figure}
\input{B4figure2.pstex_t}
\caption{}
\label{B4figure2}
\end{figure}

\begin{example}
Our final surface example is as depicted on the left-hand diagram of
Figure \ref{DegenCubic};
the singularities along the edges are of the same type as the previous
example. We take $\varphi$ to have value $0$ at $X$, $Y$, and $U$,
and $1$ at $Z$.
There are $9$ elements of the structure: six of these
are the usual ones emanating from the singularities (the dotted lines in the 
figure indicates those rays emanating from the singularities which pass into
the interiors of two-cells), and three additional
rays emanate from $U$ in the directions of $X$, $Y$, and $Z$, to produce
a consistent structure. For example, there is a ray stretching from $U$ to 
$X$ with attached function $1+yz$.
See \cite{GSInv}, Examples 4.3 and 4.4, for the
details. We would like to determine an equation for this
family in $\PP^3$ by computing the product $XYZ$. Keep in mind that we
are using the notation $X,Y,$ etc.\ for $\vartheta_X^{[1]}$, $\vartheta_Y^{[1]}$
etc.  We do this in two steps.
First we compute, say, $XY$. There will be two level two theta functions
contributing to this product, which we shall write as 
$\vartheta^{[2]}_{(X+Y)/2}$,
the theta function given by the mid-point of the line segment $XY$, and
$\vartheta^{[2]}_{U}$, the degree two theta function corresponding to the
point $U$, thought of as a half-integral point. In this notation, we in
fact have
\[
XY=\vartheta^{[2]}_{(X+Y)/2}+t \vartheta^{[2]}_{U}.
\]
One checks that $U^2=\vartheta^{[2]}_{U}$ (such a statement need not in
general be true).
We then compute $\vartheta^{[2]}_{(X+Y)/2}\cdot Z$, and find this has four 
terms: 
\[
(t+t^2)\vartheta^{[3]}_{U}+t(\vartheta^{[3]}_{(2U+X)/3}+
\vartheta^{[3]}_{(2U+Y)/3})=(t+t^2)U^3+tU^2(X+Y).
\]
The first term $tU^3$ is
purely toric; the second term $t^2U^3$ appears for the following reason.
In order to correctly compute the product using jagged paths, since $U$
appears on a number of rays, we need to perturb the endpoint a little bit.
If we perturb the endpoint so the jagged path from $Z$ now crosses
the edge $\overline{UX}$, we need to consider possible bending along
this line segment. The expression \eqref{bendingexpansion} arising
from crossing the rays on this line segment near $U$ is $(1+x)(1+yz)
=1+x+yz+t$. The monomial $t$ does not change the direction of the
jagged path, but produces the contribution $t^2U^3$.
The third and fourth terms arise as shown in the right-hand
diagram of
Figure \ref{DegenCubic}. Finally, $\vartheta^{[2]}_{U}\cdot Z=U^2Z$ is
purely toric. This gives the equation
\[
XYZ=t\big( (1+t)U^3+(X+Y+Z)U^2\big).
\]
\end{example}

\begin{figure}
\input{DegenCubic.pstex_t}
\caption{}
\label{DegenCubic}
\end{figure}

\begin{figure}
\input{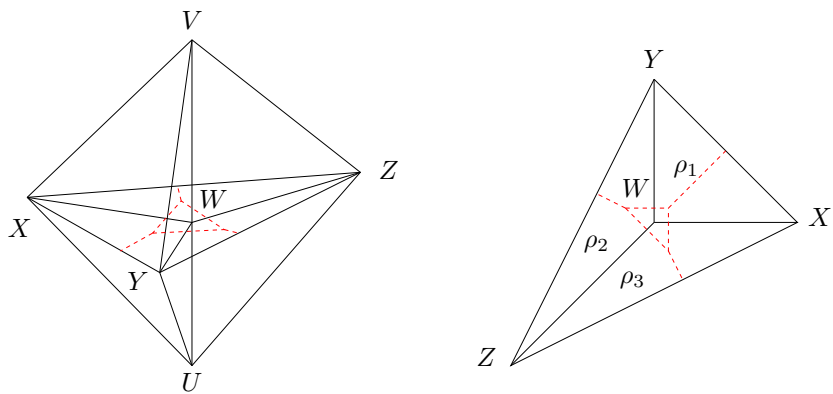}
\caption{The tropical manifold appears on the left, with the central
triangle containing the discriminant locus appearing on the right. These
three cells comprise the slabs.}
\label{MirrorLocalP2}
\end{figure}

\begin{example}
\label{crucial3dexample}
We consider one crucial three-dimensional example, see Example 5.2
of \cite{GSInv}. See Figure \ref{MirrorLocalP2}. 
As explained in \cite{GSInv}, the structure now consists of
three slabs, the two-dimensional cells depicted on the right in
the figure. The function attached to the slabs depends on which connected
component of the complement of the discriminant locus we are on. The crucial
point is the function attached to the central component, which can be written
as 
\[
f=1+x+y+z+g(xyz).
\]
Here $x,y$ and $z$ are the monomials corresponding to primitive tangent vectors
pointing from $W$ to the points $X,Y$ and $Z$ respectively. The formal
power series $g$ is determined by a procedure called \emph{normalization}
in \cite{GS11}. It must be chosen so that $\log f$ is free of pure
powers of $t=xyz$. This is expressed in terms of power series, i.e.,
that
\[
\sum_{k\ge 1} {(-1)^{k+1}\over k} (x+y+z+g(xyz))^k \in \CC\lfor x,y,z\rfor
\]
does not contain any monomial $(xyz)^l=t^l$. This determines $g$
uniquely, and one can compute $g(t)$ inductively:
\[
g(t)=-2t+5t^2-32t^3+286t^4-3038t^5+\cdots.
\]
One can then use jagged paths to determine products. The product $XYZ$
is purely toric, giving $tW^3$, while the product $UV$ can then be
written as
\[
UV=t^2(X+Y+Z+(1+g(t))W)W.
\]
Each of these terms correspond to a different term in $f$. 

Those terms
of the form $g_nt^n W^2$ will correspond to two jagged paths which
do not bend, the line segments $\overline{VW}$ and $\overline{UW}$.
As we know, these two jagged paths 
should correspond to a holomorphic triangle in the
mirror manifold. As explained in \cite{GSInv}, the mirror manifold
is an open subset $X$ of the total space of the canonical bundle of $\PP^2$.
In fact $X$ contains the zero-section of the canonical bundle, isomorphic
to $\PP^2$. The holomorphic triangle in question can be seen to intersect
the $\PP^2$ in one point, say $x$. How then do we explain the adjustments coming
from the terms of $g$? The point is that $\PP^2$ contains many holomorphic
rational curves of degree $n$ passing through $x$; by gluing any one of these
curves to the holomorphic triangle, we get a (degenerate) triangle which
should also contribute to the Floer product of $U$ and $V$. This led us 
to conjecture that the coefficient of $t^d$ in $g$ should represent a type
of $1$-point invariant for rational curves of degree $d$ in $\PP^2$.
This conjecture was supported by the observation that the sequence
of numbers $-2,5,-32,\ldots$ already appeared in several places in the
literature. First, they appeared as Gromov-Witten invariants for 
certain curve classes 
on the total space of the canonical bundle of the blow up of $\PP^2$
as given in \cite{CKYZ}. Second, they appeared in Table 6 of \cite{AKV} as
open Gromov-Witten invariants for the total space of the
anti-canonical bundle of $\PP^2$, and these numbers arose precisely
from relative homology classes of holomorphic disks where the holomorphic
disks were likely to be represented by a single disk meeting $\PP^2
\subset K_{\PP^2}$ along with a sphere contained in the $\PP^2$ attached
to the disk. 

More recently, this was argued from a different point of view, that
of Landau-Ginzburg potentials, in \cite{CLL}. This conjecture
has now been proved by Chan, Lau and Tseng in \cite{CLT}.
\end{example}

\subsection{Broken lines and functions on $\shL^{-1}$}

A jagged path is designed to organize the propagation of local
monomial sections of $\shL^{\otimes\ell}$ for any $\ell\neq 0$. Indeed, a jagged
path carries a section of $\widetilde\shP$, which in turn defines a
monomial section of $\shL^{\otimes\ell}$ in a local chart. At a bend the changes
to the section of $\widetilde\shP$ lie in $\shP$, the kernel of the
degree homomorphism $\widetilde\shP\to \ul\ZZ$. In a local chart for
$\X$ this means multiplication of the monomial section by a
monomial. Now specializing to the case $\ell=0$ we arrive at the notion
of \emph{broken lines}, which control the propagation of monomials
on $\X$. Broken lines were introduced in the literature before jagged paths
(\cite{PP2mirror},\cite{CPS},\cite{GHKI}), where they have been
used notably in the construction of Landau-Ginzburg potentials. However,
jagged paths first appeared in discussions between the two authors
and Mohammed Abouzaid in 2007. For
the following definition we use the notation from
Definition~\ref{jaggeddef}. We now also admit unbounded jagged
paths, but still with only finitely many bends. The domain of
definition of $\gamma$ is then an interval $I\subsetneq\RR$ rather than
$[0,1]$. Note that in the unbounded case there is a unique maximal
unbounded domain of linearity $(-\infty, t_0)$. 

\begin{definition}
A (bounded or unbounded) jagged path $\big(\gamma:I\to B,(m_L) \big)$
in $B$ with respect to the structure $\foD$ is called a \emph{broken
line} if for one (hence for every) domain of linearity
$L\subset I$ it holds $\deg(m_L)=0$.
\end{definition}

Broken lines are conceptually somewhat easier since it suffices to
work with the sheaf $\shP$ rather than with $\shP$ and
$\widetilde\shP$. In particular, this removes one layer of
notation.  Note also that $\textbf{vect}(m_L)$ lies in $\Lambda$ and
stays constant along a domain of linearity. Since
$\textbf{vect}(m_L)$ is negatively proportional to $\gamma'(t)$ this
shows that unlike jagged paths, broken lines can travel only in
rational directions.  

Now a little trick allows one to completely replace jagged paths by
broken lines on the technical level, and this is what is done at
most places in \cite{GHKS}. The trick is based on the observation
that a section of $\shL^{\otimes\ell}$ over an open set $U$ is the same as a
regular function over the preimage of $U$ on the total space
$\operatorname{Tot}(\shL^{-1})$ that is fibrewise homogeneous of
degree $d$.

The point is that $\operatorname{Tot}(\shL^{-1})$ has a simple
realization in terms of our program. Let $(B,\P)$ be an integral affine
manifold with singularities with polyhedral decomposition $\P$.

\begin{definition}
The \emph{truncated cone over $B$} is the integral affine manifold
defined as a set by
\[
\trcone B:= B\times [1,\infty),
\]
endowed with the following affine structure. For $\psi:U\to \RR^n$
an affine chart for $B$ defined on an open set $U\subset B$, we
define the chart
\begin{equation}\label{chart for CB}
\tilde\psi: \trcone U\lra \RR^{n+1},\quad
(x,h)\longmapsto (h\cdot \psi(x),h)
\end{equation}
for $\trcone B$. The polyhedral decomposition $\trcone\P$ is given
by the cells $\trcone\sigma:=\sigma\times[1,\infty)$ for
$\sigma\in\P$.
\end{definition}

While $\trcone B$ is topologically a product, the affine structure
is that of a cone over $B$ with tip chopped off, see
Figure~\ref{fig: truncated cone}. Note that $\trcone B$ has boundary
$\trcone(\partial B)\cup B\times\{1\}$. \begin{figure}
\input{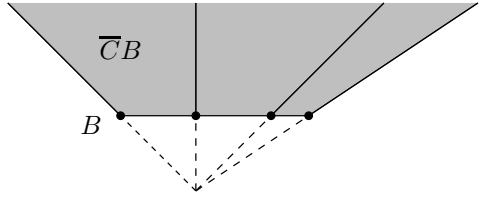} \caption{The truncated cone.} \label{fig:
truncated cone} \end{figure} As a manifestation of the cone
structure let us look at parallel transport of the vertical tangent
vector $(0,d)\in T_{(x,h)} \trcone B= T_x B\oplus\RR$ along a
straight line to $(y,h)\in \trcone B$. In a chart~\eqref{chart for
CB} centered at $x$ this tangent vector maps to $(0,d)$, but at $y$
the preimage of $(0,d)$ is $\big(-d v/h,d)$ where $v=\psi(y)-\psi(x)$.
For $h=1$, this is just the parallel transport underlying the sheaf
$\shAff(B,\ZZ)^*$! To discuss the relation with $B$ let us use
$\trcone B$ as index for $\shP$ and $\Lambda$ to distinguish these
sheaves from the corresponding sheaves on $B$.  The discussion
of parallel transport shows that $\Lambda_{\trcone B}$ restricted to
$B\times\{1\}\subset \trcone(\partial B)$ is canonically isomorphic to
the sheaf $\shAff(B,\ZZ)^*$ on $B$. Similarly,  $\widetilde\shP$,
viewed as a sheaf on $B\times\{1\}$, extends to the sheaf
$\shP_{\trcone B}$.

Diagram~\eqref{bigcommdiag} also has a simple interpretation in terms
of $\trcone B$. While there is no natural affine map from $\trcone
B$ to $B$, the projection to the $\RR$-factor defines an affine map
$\trcone B\to [1,\infty)$. This induces a homomorphism
$\Lambda_{\trcone B}\to\ul \ZZ$ and, by composition with
$\shP_{\trcone B}\to \Lambda_{\trcone B}$, a homomorphism
$\shP_{\trcone B}\to \ul\ZZ$. Thus the second column of
Diagram~\eqref{bigcommdiag} on $B$ is just the restriction of the first
column on $\trcone B$ to $B\times\{1\}$. The lower two rows are
defined by projection to $[1,\infty)$. The kernel of this projection
defines horizontal elements in $\shP_{\trcone B}$ and on
$\Lambda_{\trcone B}$. The left column is thus just the restriction of
the middle column to horizontal elements.

It is then immediate that there is a one-to-one correspondence
between broken lines on $\trcone B$ and jagged paths on $B$, simply
by composing with the projection $\trcone B\to B$ along the lines
emanating from the tip of the cone. The geometric
interpretation of this correspondence is very transparent on the
complex side.

\begin{proposition}
If $\pi:\check \X\to\Spec\CC\lfor t\rfor$ is a deformation
associated to $(B,\P)$ from Theorem~\ref{mainsmoothingtheorem} and
$\shL$ the relatively ample line bundle, then the total space of
$\shL^{-1}$ can be constructed by applying our construction to
$(\trcone B,\trcone \P)$.
\end{proposition}

To define a regular function via broken lines one needs to look at 
the asymptotic integral affine manifold of a non-compact affine
manifold, defined by equivalence classes of affine rays in unbounded
cells, the equivalence generated by affine translations, see \cite{CPS}. 
In the case
of $\trcone B$ the asymptotic integral affine manifold is just $B$.
Then given $\ell\in\NN\setminus\{0\}$ and a $1/\ell$-integral point $m$ in
the asymptotic affine manifold, consider broken lines $\big(\gamma:
(-\infty,0]\to\trcone B, (m_L)\big)$ whose unbounded asymptotic
direction is $m$, and for
$L$ the unbounded domain of linearity, $m_L$ is a $\ell$-fold multiple
of an primitive element of vanishing $t$-order. Note this fixes
$m_L$ uniquely. Now the same procedure as in the construction of
$\vartheta_m$ (Definition~\ref{Def: Lift(m)}), but with jagged paths
from $m$ to $x$ replaced by broken lines with the asymptotics
defined by $m$ and ending at $x$, defines a regular function on the
total space of the associated degeneration of non-compact varieties.
In the situation of $\trcone B$ this defines $\vartheta_m$, viewed as a
regular function on $\operatorname{Tot}(\shL^{-1})$. Note that this
regular function is fibrewise homogeneous of degree $\ell$ since for
any domain of linearity $L$, the projection of $m_L$ to $[1,\infty)$
has length $\ell$.

\section{Tropical Morse trees}

We have seen that jagged paths can be used to compute products of theta
functions, emulating the use of holomorphic triangles to compute the
product in Floer homology. This raises the question as to whether
one can compute the higher $A_{\infty}$ operations using a similar
strategy. This is explored in work in progress of Abouzaid, Gross, and
Siebert \cite{AGS}. This idea was already explored in \cite{DBr}, Chapter 8
in the case of elliptic curves (where $B=\RR/d\ZZ$ for some positive
integer $d$). In that case, as there are no singularities, one could just
use trees composed of straight lines; in the general case, one needs
to use jagged paths. We outline this here.

Let us return to the situation of \S1, where given an integral affine
manifold with singularities $B$, we assume we have $X(B)\rightarrow B$
and sections $L_\ell$ (the image of $\sigma_{\ell}$) which hopefully become
Lagrangian after a suitable choice of symplectic structure on $X(B)$. 
Further, we expect $L_\ell$ to be mirror to $\shL^{\otimes \ell}$ on 
$\check\shX$. 

Suppose we wish to compute
\begin{equation}
\label{mud2}
\mu_d:\Hom(L_{\ell_{d-1}},L_{\ell_d})\otimes\cdots\otimes\Hom(L_{\ell_0},\
L_{\ell_1})
\rightarrow\Hom(L_{\ell_0},L_{\ell_d})[2-d].
\end{equation}
As this should be defined using holomorphic disks, we hope to be able
to compute this using jagged paths instead, using the philosophy
that jagged paths correspond to holomorphic
objects which can be glued together. This leads us to the definition of
\emph{tropical Morse tree}.

To give this definition, assume given $B$ and a consistent structure
$\foD$, so that we can talk about jagged paths on $B$ with respect
to $\foD$. Further, recall that a \emph{ribbon tree} is a tree $S$ with a cyclic
ordering of edges adjacent to each vertex. This provides a cyclic ordering
of leaves of the tree, and by choosing one leaf as an output, labelled
$v_{0,d}$, and orienting each edge towards this output, 
we obtain a directed tree and canonical labelling of all other
leaves as $v_{0,1},\ldots,v_{d-1,d}$, see Figure \ref{ribbontree}.
In addition, given a sequence of integers $\ell_0,\ldots,\ell_d$, we can assign
an integer $\ell_e$ to each edge $e$ of $S$: for $e$ adjacent to $v_{i,i+1}$,
$\ell_e=\ell_{i+1}-\ell_i$. If $e$ is the outgoing edge at an interior vertex
with incoming edges $e_1,\ldots,e_p$, then $\ell_e=\sum_{i=1}^p \ell_{e_i}$.

\begin{figure}
\input{ribbontree.pstex_t}
\caption{}
\label{ribbontree}
\end{figure}

\begin{definition}
\label{TMTdef}
Suppose given distinct integers $\ell_0,\ldots,\ell_d$. Then
a tropical Morse tree with respect to the data $B,\foD$ and $\ell_0,\ldots,
\ell_d$
is a map $\psi:S\rightarrow B$ with $S$ a ribbon tree with $d+1$ leaves
whose restriction to any edge is a jagged path
(and hence comes with the additional data of attached monomials). This data
should satisfy the following conditions:
\begin{enumerate}
\item 
\[
\psi(v_{i,i+1})=p_{i,i+1}\in B\left({1\over \ell_{i+1}-\ell_i}\ZZ\right)
\]
and the initial monomial attached to the edge $e_{i,i+1}$
adjacent to $v_{i,i+1}$, viewing $\psi|_{e_{i,i+1}}$ as a jagged path, is
$z^{(p_{i,i+1})_{\varphi}}$.
\item Let $v$ be an internal vertex of $S$ with incoming edges $e_1,
\ldots,e_p$ and outgoing edge $e_{\out}$. Let $c_1z^{m_1},\ldots,c_pz^{m_p}$
be the monomials attached to the last linear segment of each jagged
path $\psi|_{e_1},\ldots,\psi|_{e_p}$. Then the monomial $c_{\out}z^{m_{\out}}$
attached to the initial linear segment of $\psi|_{e_{\out}}$ is
$\prod_{i=1}^p c_iz^{m_i}$.
\item Let $e_{0,d}$ be the edge adjacent to $v_{0,d}$, and let $c_{0,d}
z^{m_{0,d}}$ be the monomial attached to the last linear segment of
$\psi|_{e_{0,d}}$ as a jagged path. Then at $\psi(v_{0,d})$, we have
${\bf vect}(\tilde r(m_{0,d}))=0$.
\end{enumerate}
\end{definition}

Let us note a number of features of this definition. First, morally such
a tropical Morse tree should contribute to $\mu_d$ as in \eqref{mud2}
where the inputs are intersections points of Lagrangian sections
corresponding to  $p_{0,1},\ldots,p_{d-1,d}$.
However, because \eqref{mud2} is a map of degree $2-d$, this will always
be zero if $d>2$ and all inputs are of degree $0$. However, we will allow
$\ell_{i+1}-\ell_i$ to be negative, so that $\Hom(L_{\ell_i},L_{\ell_{i+1}})
\cong H^*(\check\shX,\shL^{\otimes (\ell_{i+1}-\ell_i)})$ consists only of 
top degree cohomology. Indeed, as $\shL$ is ample, it follows by
Kodaira vanishing and Serre duality that for $d<0$
\[
H^i(\shX,\shL^{\otimes \ell}) \cong {}  H^{\dim B-i}(\check\shX,
\shL^{\otimes (-\ell)})
=  \begin{cases}
H^0(\check\shX,\shL^{\otimes (-\ell)})& i=\dim B\\
0 & i<\dim B
\end{cases}
\]
Thus we can write $\{\vartheta_p \,|\, p\in B({1\over \ell}\ZZ)\}$
as a basis of $H^{\dim B}(\shX,\shL^{\otimes \ell})$ Serre dual
to the basis $\{\vartheta_p\,|\,p\in B({1\over -\ell}\ZZ)=B({1\over \ell}\ZZ)\}$
of theta functions for $H^0(\shX,\shL^{\otimes (-\ell)})$.
This allows us to treat negative and positive powers of $\shL$, and so 
we will get non-trivial possibilities for $\mu_d$ for many different
$d$.

The next point to observe is that if $\ell_{i+1}-\ell_i<0$, then in fact 
$\psi|_{e_{i,i+1}}$ needs to be viewed as a trivial jagged path which is just
a point rather than a line segment. Indeed, in this case, near $p_{i,i+1}$,
${\bf  vect}(\tilde r((p_{i,i+1})_{\varphi}))$ points
\emph{away} from $p_{i,i+1}$. But property (1) of Definition \ref{jaggeddef}
says that this vector must be negatively proportional to $\psi'$, so it is
impossible for such a jagged path to move away from $p_{i,i+1}$.

A similar argument show that if $\ell_d-\ell_0>0$, then 
$\psi|_{e_{0,d}}$ must also
be a trivial jagged path, because ${\bf vect}(\tilde r(\cdot))$ gets
bigger, not smaller, along jagged paths with attached monomials being of
positive degree. Thus in this case we cannot achieve condition (3) 
of Definition \ref{TMTdef} unless $\psi|_{e_{\out}}$ is contracted.

Note that the condition that ${\bf vect}(\tilde r(m_{0,d}))=0$ implies
that $\tilde r(m_{0,d})$, thought of as a local section of $\shAff(B,\ZZ)^*$,
is given by $(\ell_d-\ell_0)\cdot \ev_{p_{0,d}}$. Thus $m_{0,d}-(p_{0,d})_{\varphi}$
lies in $\ker(\tilde r)\cong\ZZ$; denote this difference by $\ord(\psi)$.

Finally, we note that condition (2) of Definition \ref{TMTdef} imposes a kind
of balancing condition at the vertices, namely,
\[
{\bf vect}(\tilde r(m_{\out}))=\sum_{i=1}^p {\bf vect}(\tilde r(m_i)).
\]
So the vectors on the incoming edges determine the tangent direction of the
outgoing edge.

In theory, we would like to define
\[
\mu_d(\vartheta_{p_{d-1,d}},\ldots,\vartheta_{p_{0,1}})
=\sum_{\psi} c_{0,d}t^{\ord(\psi)}\vartheta_{p_{0,d}},
\]
where the sum is over all tropical Morse trees with $\psi(v_{i,i+1})=
p_{i,i+1}$ and $p_{0,d}$ defined to be $\psi(v_{0,d})$.
We note the conditions imply that $p_{0,d}\in B\left({1\over \ell_d-\ell_0}\ZZ
\right)$. 

For $d=2$, this formula in fact recovers the theta multiplication formula
of \eqref{thetamultformula2}. Indeed, in this case $S$ is just a trivalent
tree with $3$ leaves and one vertex; the outgoing edge is necessarily 
contracted, and the balancing condition ${\bf v}_1(1)+{\bf v}_2(1)=0$
is enforced by (2) and (3) of Definition \ref{TMTdef}. 

There is a problem with this definition for $d>2$, however. The moduli spaces
of tropical Morse trees need not be the correct dimension. This happens,
for example, even in the case that $B=\RR^n/\Gamma$ is a torus, and all
the input points are contained inside a hyperplane in $\RR^n$. As pointed out
to us by M.\ Slawinski, this is already a problem when $n=1$, the case
of an elliptic curve, if all inputs are the same point in $B$. (So in fact
the arguments in \cite{DBr}, Chap.\ 8, are not complete.)

As a consequence, in order to properly define $\mu_d$ for $d\ge 3$, one
needs to perturb the moduli problem so we get finite counts of trees
when needed. We hope to find a way of doing this which preserves the
combinatorial nature of the construction, allowing for actual computations
of $A_{\infty}$-structures. This will certainly involve choices, but the
resulting $A_{\infty}$-structures should be impervious to these choices,
up to quasi-isomorphism.

Once such a choice is made, it should be easy to show that the 
$A_{\infty}$-precategory whose objects are powers of $\shL$ and whose morphisms
are given, say, by \v Cech complexes computing cohomology of powers
of $\shL$ (see \cite{DBr}, \S8.4.5 for the elliptic curve case) is 
quasi-isomorphic to an $A_{\infty}$-precategory defined using the $\mu_d$'s
along the lines of \cite{DBr}. A much more challenging problem is to
then relate this category to the actual Fukaya category of the mirror.

We also note that M.\ Slawinski, in his thesis \cite{Sl12}, introduced the
notion of tropical Morse graph with an aim of identifying a quantum
$A_{\infty}$-category structure \cite{B07} on the category
of powers of $\shL$.

\begin{examples}
(1) Consider the tree in Figure \ref{TMTfigure1}
with either $B=\RR$ or $B=\RR/m\ZZ$ for some positive
integer $m$, the latter by considering the lift of $\psi$ to the
universal cover. This tree contributes to the coefficient of
$p_{0,3}$ in $\mu_3(p_{2,3},p_{1,2},p_{0,1})$. Here we take
$\ell_0=0$, $\ell_1=1$, $\ell_2=3$ and $\ell_3=2$, noting $e_{2,3}$ and $e_{0,3}$
are contracted to points.

\begin{figure}
\input{TMTfigure1.pstex_t}
\caption{}
\label{TMTfigure1}
\end{figure}

(2) In Figure \ref{TMTfigure2}, we give a two-dimensional example,
again in $B=\RR^2$ or $\RR^2/\Gamma$ for a lattice $\Gamma$, contributing
to $\mu_4$. In this example we take $\ell_0=0$, $\ell_1=4$, $\ell_2=-4$,
$\ell_3=-7$, $\ell_4=-3$. Again, $e_{1,2}$ and $e_{2,3}$ are contracted.
\qed
\end{examples}

\begin{figure}
\input{TMTfigure2.pstex_t}
\caption{}
\label{TMTfigure2}
\end{figure}

\section{Applications of theta functions to mirror constructions}

Theta functions play a crucial role in extensive new work, partly
of the first author 
jointly with Hacking and Keel, \cite{GHKI},\cite{GHKII} 
and partly of both authors
again jointly with Hacking and Keel \cite{GHKS}. We will only discuss this work
briefly here, and only one aspect of this work.

In particular, we explain how theta functions allow us to greatly expand
the class of singularities of affine manifolds we can treat. In this survey,
we have almost exclusively considered only one two-dimensional singularity
known from the integrable systems literature as a \emph{focus-focus
singularity}. This singularity has the feature that the monodromy in
the local system $\Lambda$ about the singularity is $\begin{pmatrix}
1&1\\0&1\end{pmatrix}$, and the invariant tangent direction is tangent to
a line passing through the singularity. Further, the singularity must appear
in the interior of an edge, and not at a vertex. In higher dimensions,
the allowable singularities are the same generically: Example 
\ref{crucial3dexample} is a typical three-dimensional example, with
trivalent vertices for the discriminant locus. Theorem 
\ref{mainsmoothingtheorem} as proved in \cite{GS11} holds for this
sufficiently nice class (and a little more generally), 
called ``simple'' in \cite{GS06}. This class of
singularities is related to the structure of the degenerations
constructed by Theorem \ref{mainsmoothingtheorem}. These degenerations
are special kinds of degenerations we call \emph{toric degenerations}. These
are degenerations of Calabi-Yau varieties $f:\shX\rightarrow D$ such
that the central fibre $\shX_0$ is a union of toric varieties glued
along toric strata, and $f$ is given locally in a neighbourhood of
the most singular points of $\shX_0$ by a monomial in a toric variety.

This is an ideal class of degenerations for studying mirror symmetry,
as it exhibits the greatest level of symmetry (under mirror symmetry,
the data of the irreducible components is exchanged with the structure
of the family at the most singular points of the central fibre). 
Furthermore, it works very
well for complete intersections in toric varieties, see e.g., \cite{Gr05}.
However, one would ideally like to construct a mirror for any maximally
unipotent degenerating family $\shX\rightarrow D$, and it might be
difficult to find a birationally equivalent family which is a toric
degeneration. Thus it is desirable to expand the class of allowable
degenerations, and this is equivalent to expanding the class of allowable
singularities that our program can handle.

Let us consider the setup of \cite{GHKI} by way of example. Consider $(Y,D)$
with $Y$ a non-singular projective rational surface and $D\in |-K_Y|$ a
cycle of rational curves. We call such data a \emph{Looijenga pair}.
We can construct an integral affine manifold homeomorphic to $\RR^2$
with one singularity associated to $(Y,D)$, as follows. Let $D=D_1+\cdots+D_n$, with
$D_1,\ldots,D_n$ cyclically ordered.

For each node
$p_{i,i+1} := D_i \cap D_{i+1}$ of $D$ we take a rank
two lattice $M_{i,i+1}$ with basis $v_i,v_{i+1}$, and
the cone $\sigma_{i,i+1} \subset M_{i,i+1}
\otimes_{\ZZ} \RR$ generated by $v_i$ and $v_{i+1}$.
We then glue $\sigma_{i,i+1}$ to $\sigma_{i-1,i}$ along the rays
$\rho_i := \RR_{\geq 0} v_i$ to obtain a piecewise-linear
manifold $B$ homeomorphic
to $\RR^2$ and a decomposition
\[
\Sigma=\{\sigma_{i,i+1}\,|\,1\le i\le n\}
\cup \{\rho_i\,|\,1\le i\le n\} \cup \{0\}.
\]
We define an integral affine
structure on $B \setminus \{0\}$. We do this by defining charts
$\psi_i:U_i\rightarrow M_{\RR}$ (where $M=\ZZ^2$). Here
\[
U_i=\Int(\sigma_{i-1,i}\cup \sigma_{i,i+1})
\]
and $\psi_i$ is defined on the closure of $U_i$ by
\[
\psi_i(v_{i-1})=(1,0),\quad \psi_i(v_i)=(0,1),\quad
\psi_i(v_{i+1})=(-1,-D_i^2),
\]
with $\psi_i$ linear on $\sigma_{i-1,i}$ and $\sigma_{i,i+1}$.
The idea behind this formula is that we are pretending that 
$(Y,D)$ is in fact a toric pair. Given a ray in a two-dimensional
fan generated by $(0,1)$ corresponding to a divisor $C$, with adjacent
rays generated by $(1,0)$ and $(-1,-D_i^2)$ respectively, one has
$C^2=-D_i^2$. In particular, if $(Y,D)$ were in fact toric,
the above construction would just yield $B\cong\RR^2$ as an affine manifold,
with $\Sigma$ the fan defining $Y$. If $(Y,D)$ is not toric,
$B$ has a non-trivial singularity at the origin.

The reader can check a simple example: if $Y$ is a del Pezzo surface of
degree $5$, one can find a cycle of $5$ $-1$-curves on $Y$ giving $D$.
In this case, the monodromy of $\Lambda$ about the resulting singularity
is $\begin{pmatrix}1&1\\ -1&0\end{pmatrix}$; see \cite{GHKI}, Example 1.8
for details.

Here $(B,\Sigma)$ can be thought of as a dual intersection complex of
$(Y,D)$. If one reinterprets $(B,\Sigma)$ as an intersection complex
for a degeneration, one would hope to find a flat family $\check\shX
\rightarrow \Spf \CC\lfor t\rfor$ whose central fibre is a union of $n$
copies of $\AA^2$. Specifically, $\check\shX_0$ should be the
\emph{$n$-vertex}, the union of coordinate planes (if $n\ge 3$)
\[
\VV_n=\AA^2_{x_1,x_2}\cup\AA^2_{x_2,x_3}\cup\cdots\cup\AA^2_{x_n,x_1}
\subseteq\AA^n,
\]
where the subscripts denote the non-zero coordinates on each plane.

One can attempt to use the techniques of \cite{GS11} to produce such
a deformation. The problem is there is no local model for a smoothing in
a neighbourhood of $0\in\VV_n$, and the arguments of \cite{GS11} work
by gluing together local models, which are required in codimension
$\le 2$. However, if we throw away $0\in\VV_n$,
a choice of a strictly convex piecewise linear function $\varphi$ on
$B$ gives rise to a $k$-th order deformation of $\VV_n^o:=\VV_n\setminus
\{0\}$, denoted $\check\shX^o_k\rightarrow \Spec\CC[t]/(t^{k+1})$.
This deformation looks purely toric in a neighbourhood of each 
connected component of
$\Sing(\VV_n^o)$. Since $\VV_n$ is affine, we can try to recover
a $k$-th order deformation $\check\shX_k$ of $\VV_n$ as follows. Suppose
the coordinates $x_1,\ldots,x_n$ on $\VV_n$
lift to functions on $\check\shX_k^o$. Then we can embed $\check\shX_k^o$
into $\AA^n\times\Spec\CC[t]/(t^{k+1})$ using these lifts, and take 
the closure of the image to be $\check\shX_k$. The problem is that for
the naive family $\check\shX^o_k$, the $x_i$ won't lift.

The solution is to modify the construction of $\check\shX_k^o$ via
a structure as in \cite{GS11}: this structure should consist of lines radiating
from the origin. The functions attached to these lines determine 
automorphisms, and these automorphisms are then used to modify the
standard gluings, giving a different deformation $\check\shX_k^o$. Any
structure will provide such a deformation, as there is no meaningful
compatibility of automorphisms which can be checked at the origin, as
we have deleted the origin. If we had not deleted it, there is
no local model for a smoothing there in which we could have checked such
compatibility.

However, one can instead insist on using a structure which is consistent,
which allows us to define theta functions on $\check\shX_k^o$, which yield
lifts of the variables $x_i$. To do so requires one to make
a careful guess for the structure. In \cite{GHKI}, we define a canonical
choice of structure motivated by \cite{GPS}: this structure encodes
certain relative Gromov-Witten invariants of the pair $(Y,D)$. Importantly,
one needs the original surface to construct the structure; it depends
on more than just the singularity. In particular, it is possible to
have two different pairs $(Y,D)$, $(Y',D')$ giving rise to the same
affine singularity but to different structures.

This allows one to construct the deformations $\check\shX_k\rightarrow
\Spec\CC[t]/(t^{k+1})$, and taking the limit, one obtains a formal
deformation $\check\shX\rightarrow\Spf\CC\lfor t\rfor$.

The following is a more precise statement of the main result of \cite{GHKI}:

\begin{theorem} 
\label{GHKImaintheorem}
Let $(Y,D)$ be a Looijenga pair and let $\sigma_P\subseteq H_2(Y,\RR)$
be a strictly convex rational polyhedral cone containing the Mori
cone (the cone of effective curves) of $Y$. Let $P=\sigma_P\cap
H_2(Y,\ZZ)$, and $\fom_P$ the monomial ideal in $\CC[P]$ generated
by $\{z^p\,|\, p\in P\setminus\{0\}\}$. Let $\widehat{\CC[P]}$ be the
formal completion of $\CC[P]$ with respect to $\fom_P$. Then there is a formal
smoothing of the $n$-vertex $\check\shX\rightarrow\Spf\widehat{\CC[P]}$
canonically associated to $(Y,D)$. This family can be viewed as the mirror
family to the pair $(Y,D)$.

If $D$ supports a divisor ample on $Y$, then we can take $\sigma_P$ to
be the Mori cone, and the mirror family extends to a family 
$\check\shX\rightarrow\Spec\CC[P]$.
\end{theorem}

\begin{example}
The cubic surface with a triangle of $-1$-curves as $D$ provides a particularly
attractive example. In this case the monodromy
of the singularity is minus the identity. The relevant structure controlling
the deformation is extremely complicated, with there being a non-trivial
ray of every possible rational slope. Nevertheless, it can be shown that
the theta function lifts of $x_1,x_2,x_3$ satisfy a very simple cubic algebraic
equation, see \cite{GHKI}, Example 5.12. Taking the closure of this family
in $\PP^3$ gives a universal family of cubic surfaces constructed by Cayley.
\end{example}

The above results have an application to a conjecture of Looijenga
\cite{L81} predating mirror symmetry, concerning the deformation theory of cusp
singularities. A cusp singularity is a normal surface singularity whose
minimal resolution has exceptional divisor a cycle of rational curves.

It had been observed that cusp singularities come in \emph{dual
pairs}. This can be explained as follows.

Let $M=\ZZ^2$.
Let $T \in \SL(M)$ be a hyperbolic matrix, i.e., having a real eigenvalue
$\lambda >1$. 
Let $w_1,w_2 \in M_{\RR}$ be eigenvectors with eigenvalues
$\lambda_1=1/\lambda$, $\lambda_2 = \lambda$,
chosen so that $w_1 \wedge w_2 > 0$ (in the standard counter-clockwise
orientation of $\RR^2$). Let $\bar C,\bar C'$ 
be the strictly convex cones spanned
by $w_1,w_2$ and $w_2,-w_1$, and let $C,C'$ be their interiors, either of
which is preserved by $T$. Let $U_C,U_{C'}$ be the corresponding tube
domains, i.e.,
\[
U_C := \{z \in M_{\CC}| \Im(z) \in C\}/M \subset M_{\CC}/M = M \otimes \Gm.
\]
$T$ acts freely and properly discontinuously on $U_C, U_{C'}$. The
holomorphic hulls of the quotients
$U_{C}/\langle T \rangle,U_{C'}/\langle T \rangle$ are normal
surface germs. These are each a cusp singularity, and they are
considered dual to each other. All cusps (and their duals) arise this way.

We then obtain a proof of Looijenga's conjecture concerning smoothability
of cusp singularities:

\begin{theorem} A germ of a cusp singularity $(X,p)$ is smoothable
if and only if there is a Looijenga pair $(Y,D)$ along with a birational
morphism $Y\rightarrow\bar Y$ contracting $D$ to the dual cusp singularity.
\end{theorem}

Looijenga proved smoothability of $(X,p)$ implies the existence of $(Y,D)$.
This is done by realising a cusp singularity and its dual inside an
Inoue surface, and then deforming the Inoue surface so that $p$ smooths
but the dual cusp remains untouched. The resulting surface resolves
to $(Y,D)$. The converse is proved in \cite{GHKI}, by studying the mirror
family to $(Y,D)$. One can show with some additional effort 
that the family extends (analytically) to one which contains as a fibre
the dual cusp to $p$. Since the mirror family is a smoothing, the dual
cusp is thus smoothable.

There are also connections between this construction and cluster
varieties associated to skew symmetric matrices of rank $2$ (see
\cite{FG06} for definitions). In particular, theta functions suggest 
a vast generalisation of the (conjectural) Fock-Goncharov dual bases for
cluster varieties. In particular, the above construction leads
to a proof of the Fock-Goncharov conjecture for the $X$-cluster
variety associated to such a rank $2$ skew-symmetric matrix.

The ideas used to prove Theorem \ref{GHKImaintheorem} can currently
be generalized to the K3 case. One starts with a one-parameter
maximally unipotent degeneration $\shY\rightarrow D$ of K3 surfaces.
We assume $\shY$ is non-singular and the map to $D$ is normal
crossings and relatively minimal. We can then build a dual intersection complex
$(B,\P)$ of this degeneration, where $\P$ is a decomposition of
$B$ into standard simplices. All the singularities of $B$ now lie at 
vertices, reflecting the geometry of the irreducible components of $\shY_0$.
The mirror family is constructed by first building a union $\check \shX_0$
of projective planes whose intersection complex is $(B,\P)$. We smooth by
constructing a suitable structure on $B$. This will define deformations
of $\check\shX^o_0$ obtained from $\check\shX$ by deleting zero-dimensional
strata. The correct choice of structure will be consistent, yielding
well-defined theta functions on these deformations, enabling us
to compactify the deformations much as before. As a consequence, we obtain
theta functions on K3 surfaces which enjoy many nice properties; see
\cite{GHKK3} for details. Crucially, we can show that theta functions
are essentially independent of the choice of birational model
of $\shY\rightarrow D$. This leads to a proof of a strong form of Tyurin's
conjecture.

This general construction of mirrors for K3 surfaces gives encouragement
that a similar construction will work in all dimensions. At the moment,
the techniques available rely heavily on \cite{GPS}, which is a two-dimensional
result. However, we anticipate that a suitable understanding of log
Gromov-Witten invariants \cite{GSlog},\cite{Ch},\cite{ACh}  
will allow us to create consistent
structures in general. Certain types of Gromov-Witten invariants 
associated to the degeneration $\shY\rightarrow D$ will be used
to construct the structure defining the mirror. Morally, these Gromov-Witten
invariants will count holomorphic disks with boundary contained in fibres
of the SYZ fibration, but we expect a purely algebro-geometric description
of these invariants.

\end{document}